\documentclass[11pt]{amsart}
\usepackage{amsfonts}
\usepackage{mathrsfs}
\usepackage{amsmath}
\usepackage{amssymb, amsmath}
\usepackage{graphicx}

\newcommand{\newcom}{\newcommand}

\newcom{\al}{\alpha}
\newcom{\be}{\beta}
\newcom{\eps}{\epsilon}
\newcom{\veps}{\varepsilon}
\newcom{\ga}{\gamma}
\newcom{\Ga}{\Gamma}
\newcom{\ka}{\kappa}
\newcom{\Lam}{\Lambda}
\newcom{\lam}{\lambda}
\newcom{\Om}{\Omega}
\newcom{\om}{\omega}
\newcom{\Si}{\Sigma}
\newcom{\si}{\sigma}
\newcom{\tht}{\theta}
\newcom{\dtri}{\nabla}
\newcom{\tri}{\triangle}
\newcom{\oo}{\infty}
\newcom{\vphi}{\varphi}
\newcom{\cB}{{\mathcal B}}
\newcom{\cC}{{\mathcal C}}
\newcom{\cD}{{\mathcal D}}
\newcom{\cF}{{\mathcal F}}
\newcom{\cL}{{\mathcal L}}
\newcom{\cM}{{\mathcal M}}
\newcom{\cP}{{\mathcal P}}
\newcom{\cS}{{\mathcal S}}
\newcom{\cQ}{{\mathcal Q}}
\newcom{\cT}{{\mathcal T}}
\newcom{\cY}{{\mathcal Y}}
\newcom{\cZ}{{\mathcal Z}}
\newcom{\R}{\mathbb R}
\newcom{\T}{\mathbb T}
\newcom{\N}{\mathbb N}
\newcom{\Z}{\mathbb Z}
\newcom{\C}{\mathbb C}
\newcom{\E}{\mathbb E}

\newcom{\f}{\frac}
\newcom{\di}{\displaystyle\int}
\newcom{\ds}{\displaystyle\sum}
\newcom{\dl}{\displaystyle\lim}
\newcom{\ov}{\overline}
\newcom{\sset}{\subset}
\newcom{\wt}{\widetilde}
\newcom{\pa}{\partial}
\newcom{\p}{\partial}
\newcom\na{\nabla}
\newcom{\suml}{\sum\limits}
\newcom{\supl}{\sup\limits}
\newcom{\intl}{\int\limits}
\newcom{\infl}{\inf\limits}
\newcom{\disp}{\displaystyle}
\newcom{\non}{\nonumber}
\newcom{\no}{\noindent}
\newcom{\QED}{$\square$}

\newtheorem{athm}{\bf \t}[section]
\newenvironment{thm} [1] {\def\t{#1}\begin{athm} \bf \rm} {\end{athm}}
\newcom{\bthm}{\begin{thm}}
\newcom{\ethm}{\end{thm}}

\newtheorem{theorem}{Theorem}[section]
\newtheorem{lemma}{Lemma}[section]

\newtheorem{proposition}{Proposition}[section]

\newcom{\beq}{\begin{equation}}
\newcom{\eeq}{\end{equation}}
\newcom{\ben}{\begin{eqnarray}}
\newcom{\een}{\end{eqnarray}}
\newcom{\beno}{\begin{eqnarray*}}
\newcom{\eeno}{\end{eqnarray*}}
\newcom{\bali}{\begin{aligned}}
\newcom{\eali}{\end{aligned}}

\numberwithin{equation}{section}

\begin{document}

\title[Asymptotic stability of the compressible Euler-Maxwell equations ]
 {Asymptotic stability of stationary
solutions to the compressible Euler-Maxwell equations }

\author{Qingqing Liu}
\address{(QQL)
The Hubei Key Laboratory of Mathematical Physics, School of
Mathematics and Statistics, Central China Normal University, Wuhan,
430079, P. R. China} \email{shuxueliuqingqing@126.com}

\author{Changjiang Zhu*}
\address{(CJZ)
The Hubei Key Laboratory of Mathematical Physics, School of
Mathematics and Statistics, Central China Normal University, Wuhan,
430079, P. R. China} \email{cjzhu@mail.ccnu.edu.cn}
\thanks{*Corresponding author. Email: cjzhu@mail.ccnu.edu.cn }


\date{\today}


\begin{abstract}
In this paper, we are concerned with the compressible Euler-Maxwell
equations with a nonconstant background density (e.g. of ions) in
three dimensional space. There exist stationary solutions when the
background density is a small perturbation of a positive constant
state. We first show the asymptotic stability of solutions to the
Cauchy problem near the stationary state provided that the initial
perturbation is sufficiently small. Moreover the convergence rates
are obtained by combining the $L^p$-$L^q$ estimates for the
linearized equations with time-weighted estimate.
\end{abstract}

\maketitle

\tableofcontents

\section{Introduction}

The dynamics of two separate compressible fluids of ions and
electrons interacting with their self-consistent electromagnetic
field in plasma physics can be described by the compressible 2-fluid
Euler-Maxwell equations \cite{Besse,Rishbeth}. In this paper, we
consider the following one-fluid compressible Euler-Maxwell system
when the background density $n_{b}$ is a function of spatial
variable and the electrons flow is isentropic (see
\cite{Duan,UK,USK} when $n_{b}=const.$), taking the form of
\begin{eqnarray}\label{1.1}
&&\left\{\begin{aligned}
&\partial_t n+\nabla\cdot(nu)=0,\\
&\partial_t u+u \cdot \nabla u+\frac{1}{n}\nabla
p(n)=-(E+u\times B)-\nu u,\\
&\partial_t E-\nabla\times B=nu,\\
&\partial_t B+\nabla \times E=0,\\
&\nabla \cdot E=n_{b}(x)-n, \ \  \nabla \cdot B=0.
\end{aligned}\right.
\end{eqnarray}
Here, $n=n(t,x)\geq 0 $ is the  electron density, $ u=u(t,x)\in
\mathbb{R}^{3}$ is the electron velocity, $ E=E(t,x)\in
\mathbb{R}^{3}$, $ B=B(t,x)\in \mathbb{R}^{3}$, for $ t>0, \ x \in
\mathbb{R}^{3} $, denote electron and  magnetic fields respectively.
Initial data is given as
\begin{eqnarray}\label{1.2}
[n,u,E,B]|_{t=0}=[n_{0},u_{ 0},E_{0},B_0],\ \ \ x\in\mathbb{R}^{3},
\end{eqnarray}
with the compatible conditions
\begin{eqnarray}\label{1.3}
\nabla \cdot E_0=n_{b}(x)-n_{0}, \ \  \nabla \cdot B_0=0, \ \ \
x\in\mathbb{R}^{3}.
\end{eqnarray}
The pressure function $ p(\cdot)$ of the flow depending only on the
density satisfies the power law $p(n)=A n^{\gamma}$ with constants
$A>0$ and the adiabatic exponent  $\gamma >1 $. Constant $\nu>0$ is
the velocity relaxation frequency. In this paper, we set $ A=1,\
\nu=1$ without loss of generality. $n_{b}(x)$ denotes the stationary
background ion density satisfying
\begin{eqnarray*}
n_{b}(x)\rightarrow n_{\infty}, \ \ \ \  \textrm{as}\ \ \ \
|x|\rightarrow \infty,
\end{eqnarray*}
for a positive constant state $n_{\infty}>0$. Throughout this paper,
we take $n_{\infty}=1$ for simplicity.

In comparison with the Euler-Maxwell system studied in \cite{Duan},
where the background density is a uniform constant, the naturally
existing steady states of system \eqref{1.1} are no longer constants
$[1,0,0,0]$. The stationary equations to the Cauchy problem
\eqref{1.1}-\eqref{1.2} are given as
\begin{eqnarray}\label{sta.eq0}
\left\{\begin{aligned} &\frac{1}{n_{st}}\nabla
p(n_{st})=-E_{st},\\
&\nabla \times E_{st}=0,\\
&\nabla \cdot E_{st}=n_{b}(x)-n_{st}.
\end{aligned}\right.
\end{eqnarray}

 First, in this paper, we prove the existence of the stationary
 solutions to the Cauchy problem \eqref{1.1}-\eqref{1.2} under
 some conditions on the background density $n_{b}(x)$. For this
 purpose, let us define the weighted norm
 $\|\cdot\|_{W_{k}^{m,2}}$ by
\begin{eqnarray}\label{def.norm}
\|g\|_{W_{k}^{m,2}}=\left(\sum_{|\alpha|\leq
m}\int_{\mathbb{R}^{3}}(1+|x|)^{k}|\partial^{\alpha}_{x}g(x)|^2dx\right)^{\frac{1}{2}},
\end{eqnarray}
for suitable $g=g(x)$ and integers $m\geq0$, $k\geq0$.

Actually, one has the following theorem.
\begin{theorem}\label{sta.existence}
For integers $m\geq 2$ and $k\geq 0$, suppose that
$\|n_{b}-1\|_{W_{k}^{m,2}}$ is small enough. Then the stationary
problem \eqref{sta.eq0} admits a unique solution $(n_{st},E_{st})\in
L^\infty(0, T; W_{k}^{m,2})$ satisfying
\begin{eqnarray}\label{sta.pro}
\|n_{st}-1\|_{{W_{k}^{m,2}}}\leq C \|n_{b}-1\|_{W_{k}^{m,2}},\ \ \
\|E_{st}\|_{{W_{k}^{m-1,2}}}\leq C \|n_{b}-1\|_{W_{k}^{m,2}},
\end{eqnarray}
for some constant $C$.
\end{theorem}

There have been extensive investigations into the simplified
Euler-Maxwell system where all the physical parameters are set to
unity. For the one-fluid Euler-Maxwell system, by using the
fractional Godunov scheme as well as the compensated compactness
argument, Chen-Jerome-Wang in \cite{Chen} proved global existence of
weak solutions to the initial-boundary value problem in one space
dimension for arbitrarily large initial data in $L^{\infty}$. Jerome
in \cite{Jerome} established a local smooth solution theory for the
Cauchy problem over $\mathbb{R}^3$ by adapting the classical
semigroup-resolvent approach of Kato in \cite{Kato}. Recently, Duan
in \cite{Duan} proved the existence and uniqueness of global
solutions in the framework of smooth solutions with small amplitude,
moreover the detailed analysis of Green's function to the linearized
system was made to derive the optimal time-decay rates of perturbed
solutions. The similar results are independently given by
Ueda-Wang-Kawashima in \cite{USK} and Ueda-Kawashima in \cite{UK} by
using the pure time-weighted energy method. For the the original
two-fluid Euler-Maxwell systems
 with various parameters, the limits as some parameters go to zero
 have been studied recently.
 Peng-Wang in \cite{Peng ,PW1,PW2} justified the convergence of the one-fluid compressible
Euler-Maxwell system to the incompressible Euler system,
compressible Euler-Poisson system and an electron
magnetohydrodynamics system for well-prepared smooth initial data.
These asymptotic  limits are respectively called non-relativistic
limit, the quasi-neutral limit and the limit of their combination.
Recently, Hajjej and Peng in \cite{HP} considered the
zero-relaxation limits for periodic smooth solutions of
Euler-Maxwell systems. For the 2-fluid Euler-Maxwell system,
depending on the choice of physical parameters, especially the
coefficients of $u_{\pm}$ were assumed $\nu_{+}=\nu_{-}$,
Duan-Liu-Zhu in \cite{DLZ} obtained the existence and the time-decay
rates of the solutions.  Much more studies have been made for the
Euler-Poisson system when the magnetic field is absent; see
\cite{Guo,GuoPausader,Luo,Deng,Smoller,Chae} and references therein
for discussion and analysis of the different issues such as the
existence of global smooth irrotational flow \cite{Guo} for an
electron fluid and \cite{GuoPausader} for the ion dynamics, large
time behavior of solutions \cite{Luo}, stability of star solutions
\cite{Deng,Smoller} and finite time blow-up \cite{Chae}.

However, there are few results on the global existence of solutions
to the Euler-Maxwell system when the non-moving ions  provide a
nonconstant background $n_{b}(x)$, whereas in many papers related to
one-fluid Euler-Maxwell  system $n_{b}=1$. In this paper, we prove
that there exists a stationary  solution when the background density
is a small perturbation of a positive constant state and we show the
asymptotic stability of  the stationary solution and then obtain the
convergence rate of the global solution towards the stationary
solution. The main result is stated as follows. Notations will be
explained at the end of this section.

\begin{theorem}\label{Corolary}
Let $ N\geq 3$ and $ \eqref{1.3}$ hold. Suppose
$\|n_{b}-1\|_{W_{0}^{N+1,2}}$ is small enough. Then there are $
\delta_{0}>0$, $ C_{0}>0$ such that if
\begin{eqnarray*}
\|[n_{0}-n_{st},u_{0},E_{0}-E_{st},B_{0}]\|_{N} \leq \delta_{0},
\end{eqnarray*}
then, the Cauchy problem $\eqref{1.1}$-$\eqref{1.2}$ admits a unique
global solution  $[n(t,x),u(t,x),E(t,x),B(t,x)] $ satisfying
\begin{eqnarray*}
[n-n_{st},u,E-E_{st},B]\in C([0,\infty);H^{N}(\mathbb{R}^{3}))\cap
{\rm Lip}([0,\infty);H^{N-1}(\mathbb{R}^{3})),
\end{eqnarray*}
and
\begin{eqnarray*}
\sup_{t \geq 0}\|[n(t)-n_{st},u(t),E(t)-E_{st},B(t)]\|_{N}\leq C_{0}
\|[n_{0}-n_{st},u_{0},E_{0}-E_{st},B_{0}]\|_{N}.
\end{eqnarray*}
Moreover, there are $\delta_{1}>0$, $ C_{1}>0$ such that if
\begin{eqnarray*}
\|[n_{0}-n_{st},u_{0},E_{0}-E_{st},B_{0}]\|_{N+3}+\|[u_{0},E_{0}-E_{st},B_{0}]\|_{L^{1}}\leq
\delta_{1},
\end{eqnarray*}
and $\|n_{b}-1\|_{W_{0}^{N+4,2}}$ is small enough, then the solution
$[n(t,x),u(t,x),E(t,x),B(t,x)] $ satisfies that for any  $ t \geq
0$,
\begin{eqnarray}\label{UN.decay}
\|[n(t)-n_{st},u(t),B(t),E(t)-E_{st}]\|_{N} \leq C_{1}
(1+t)^{-\frac{3}{4}},
\end{eqnarray}
\begin{eqnarray}\label{UhN.decay}
\|\nabla[n(t)-n_{st},u(t),B(t),E(t)-E_{st}]\|_{N-1} \leq C_{1}
(1+t)^{-\frac{5}{4}}.
\end{eqnarray}
More precisely, if
\begin{eqnarray*}
\|[n_{0}-n_{st},u_{0},E_{0}-E_{st},B_{0}]\|_{6}+\|[u_{0},E_{0}-E_{st},B_{0}]\|_{L^{1}}\leq
\delta_{1},
\end{eqnarray*}
and $\|n_{b}-1\|_{W_{0}^{7,2}}$ is small enough, we have
\begin{eqnarray}\label{sigmau.decay}
\|[n(t)-n_{st},u(t)]\| \leq C_{1} (1+t)^{-\frac{5}{4}},
\end{eqnarray}
\begin{eqnarray}\label{EB.decay}
\|[E(t)-E_{st},B(t)]\|\leq C_{1}(1+t)^{-\frac{3}{4}}.
\end{eqnarray}
If
\begin{eqnarray*}
\|[n_{0}-n_{st},u_{0},E_{0}-E_{st},B_{0}]\|_{7}+\|[u_{0},E_{0}-E_{st},B_{0}]\|_{L^{1}}\leq
\delta_{1},
\end{eqnarray*}
and $\|n_{b}-1\|_{W_{0}^{8,2}}$ is small enough, then $E(t)$
satisfies
\begin{eqnarray}\label{E.decay}
\|E(t)-E_{st}\|\leq C_{1}(1+t)^{-\frac{5}{4}}.
\end{eqnarray}
\end{theorem}

The proof of existence in Theorem \ref{Corolary} is based on the
classical energy method. As in \cite{Duan}, the key point is to
obtain the uniform-in-time {\it a priori} estimates in the form of
$$
\mathcal{E}_N(\bar{V}(t))+\lambda
\int_0^t\mathcal{D}_N(\bar{V}(s))\,ds\leq \mathcal{E}_N(\bar{V}_0),
$$
where $\bar{V}(t)$ is the perturbation of solutions, and
$\mathcal{E}_N(\cdot)$, $\mathcal{D}_N(\cdot)$ denote the energy
functional and energy dissipation rate functional. Here if we make
the energy estimates like what Duan  did in \cite{Duan}, it is
difficult to control the highest-order derivative of $\bar{E}$
because of the regularity-loss type in the sense that
$[\bar{E},\bar{B}]$ is time-space integrable up to $N-1$ order only.
In this paper, we modify the energy estimates by choosing a weighted
function $1+\sigma_{st}+\Phi(\sigma_{st})$ which plays a vital role
in closing the energy estimates.

Furthermore,  for the convergence rates of perturbed solutions in
Theorem 1.1, we can not analyze the corresponding linearized system
of \eqref{1.1} around the steady state $[n_{st},0,E_{st},0]$
directly. In this case, the Fourier analysis fails due to the
difficulty of variant coefficients. Here, the main idea follows from
\cite{Duan} for combining energy estimates with the linearized
results in \cite{Duan}. In the process of obtaining  the fastest
decay rates of the perturbed solution, the great difficulty is to
deal with these linear nonhomogeneous sources including $\rho_{st}$,
which can not bring enough decay rates. Whereas in \cite{Duan}, the
nonhomogeneous sources are at least quadratically nonlinear. To
overcome this difficulty, we make iteration for the inequalities
\eqref{sec5.ENV0} and $\eqref{sec5.high}$ together. In theorem
\ref{Corolary}, we only capture the same time-decay properties of
$u,\ E-E_{st}$ and $B$ as \cite{Duan} except $n-n_{st}$.
$\|n-n_{st}\|$ decays as $(1+t)^{-\frac{5}{4}}$ in the fastest way,
because the nonhomogeneous sources containing $\rho_{st}$ decay at
most the same as $\sqrt{\mathcal{E}^h_N(\cdot)}$.

The similar work was done for Vlasov-Poisson-Boltzmann system, where
the background density is also a function of spatial variable. Duan
and Yang in \cite{RY} considered the stability of the stationary
states which were given by an elliptic equation with the exponential
nonlinearity. The optimal time-decay of the Vlasov-Poisson-Boltzmann
system in $\mathbb{R}^{3}$ was obtained by Duan and Strain in
\cite{DS}. We also mention the work Duan-Ukai-Yang-Zhao in
\cite{RSYZ}, Duan-Liu-Ukai-Yang in \cite{DLUY} for the study of
optimal convergence rates of the compressible Navier-Stokes
equations with potential forces. Their proofs were based on the
combination of spectral analysis and energy estimates. Recently,
Duan-Ukai-Yang in \cite{RSY} developed a method of the combination
of the spectral analysis and energy estimates to deal with the
optimal time decay for study of equations of gas motion.

We further remark the result in \cite{RY}, the existence of solution
to the elliptic equation $\Delta \phi=e^{\phi}-\bar{\rho}(x)$ has
been proved when $\|\bar{\rho}-1\|_{W_{k}^{m,\infty}}$ is
sufficiently small, where the weighted norm
$\|\cdot\|_{W_{k}^{m,\infty}}$ is defined by
\begin{eqnarray}\label{def.norm1}
\|g\|_{W_{k}^{m,\infty}}=\sup_{x\in\mathbb{R}^{3}}(1+|x|)^{k}\sum_{|\alpha|\leq
m}|\partial^{\alpha}_{x}g(x)|
\end{eqnarray}
for suitable $g=g(x)$ and integers $m\geq0$, $k\geq0$, the stability
of the perturbed solutions can be proved when
$\|\bar{\rho}-1\|_{W_{2}^{N+1,\infty}}$ is sufficiently small. We
can also prove the stability of stationary solutions in the
framework of \cite{RY} if $\|n_{b}-1\|_{W_{0}^{N+1,\infty}}$ is
sufficiently small. In order to obtain the same convergence rates,
$\|n_{b}-1\|_{W_{2}^{N+4,\infty}}$ should be sufficiently small in
the process of dealing with $\rho_{st}\bar{u}$ as in Section
\ref{sec4},
$$
\|\rho_{st} \bar{u}\|_{L^1} \leq \|\rho_{st}\|\left\|\bar{u}
\right\| \leq C \|\rho_{st}\|_{W_{2}^{N+4,\infty}}\| \bar{u}\|.
$$
Notice that $W_{2}^{N+4,\infty}\subseteq W_{0}^{N+4,2}$, it seems to
be better to consider the existence of steady states in the weighted
Sobolev space $W_{k}^{m,2}$.

 Let us introduce some notations for the use throughout this paper. $C$ denotes some
positive (generally large) constant and $ \lambda$ denotes some
positive (generally small) constant, where both $C$ and $ \lambda$
may take different values in different places. For two quantities
$a$ and $b$, $a\sim b$ means $\lambda a \leq  b \leq
\frac{1}{\lambda} a $ for a generic constant $0<\lambda<1$. For any
integer $m\geq 0$, we use $H^{m}$, $\dot{H}^{m}$ to denote the usual
Sobolev space $H^{m}(\mathbb{R}^{3})$ and the corresponding
$m$-order homogeneous Sobolev space, respectively. Set $L^{2}=H^{m}$
when $m = 0$. For simplicity, the norm of $ H^{m}$ is denoted by
$\|\cdot\|_{m} $ with $\|\cdot \|=\|\cdot\|_{0}$. We use $
\langle\cdot, \cdot \rangle$ to denote the inner product over the
Hilbert space $ L^{2}(\mathbb{R}^{3})$, i.e.
\begin{eqnarray*}
\langle f,g \rangle=\int_{\mathbb{R}^{3}} f(x)g(x)dx,\ \ \ \  f =
f(x),\ \  g = g(x)\in L^2(\mathbb{R}^{3}).
\end{eqnarray*}
 For a multi-index $\alpha =
[\alpha_1, \alpha_2, \alpha_3]$, we denote $\partial^{\alpha} =
\partial^{\alpha_{1}}_ {x_1}\partial^{\alpha_{2}}_ {x_2} \partial^{\alpha_{3}}_ {x_3} $.
The length of $ \alpha$ is $|\alpha| = \alpha_1 + \alpha_2 +
\alpha_3$. For simplicity, we also set
$\partial_{j}=\partial_{x_{j}}$ for $j = 1, 2, 3$.

We conclude this section by stating the arrangement of the rest of
this paper. In Section 2, we prove the existence of the stationary
solution. In Section 3, we reformulate the Cauchy problem under
consideration and obtain asymptotic stability of solutions near the
stationary state provided that the initial perturbation is
sufficiently small. In Section 4, we study the time-decay rates of
solutions to the stationary solutions by combining the $L^p$-$L^q$
time-decay property of the linearized homogeneous system with
time-weighted estimate.

\vspace{5mm}

\section{Existence of stationary solution}\label{sec2}
In this section, we will prove the existence of stationary solutions
to $\eqref{sta.eq0}$ by using the contraction mapping theorem. From
$\eqref{sta.eq0}_2$, there exists $\phi_{st}$ such that
$E_{st}=\nabla \phi_{st}$, it turns equation $\eqref{sta.eq0}$ into
\begin{eqnarray}\label{sta.eq1}
\left\{\begin{aligned} &\frac{1}{n_{st}}\nabla
p(n_{st})=-\nabla\phi_{st},\\
&\Delta\phi_{st}=n_{b}(x)-n_{st}.
\end{aligned}\right.
\end{eqnarray}
We introduce the nonlinear transformation (cf. \cite{Deng})
\begin{eqnarray}\label{sta.tra}
Q_{st}=\frac{\gamma}{\gamma-1}(n_{st}^{\gamma-1}-1).
\end{eqnarray}
From $\eqref{sta.eq1}$ and $\eqref{sta.tra}$, we derive the
following elliptic equation
\begin{eqnarray}\label{sta.ellip}
\Delta
Q_{st}=\left(\frac{\gamma-1}{\gamma}Q_{st}+1\right)^{\frac{1}{\gamma-1}}-n_{b}(x).
\end{eqnarray}
For convenience, we replace $Q_{st}$ by $\phi$ in the following.
Equation $\eqref{sta.ellip}$ can be rewritten as the integral form
\begin{eqnarray*}
\phi=T(\phi)=G*\left(\left(\frac{\gamma-1}{\gamma}\phi+1\right)^{\frac{1}{\gamma-1}}-\frac{1}{\gamma}\phi
-n_{b}(x)\right),
\end{eqnarray*}
where $G=G(x)$ given by
\begin{eqnarray*}
G(x)=-\frac{1}{4\pi|x|}e^{-\tfrac{1}{\sqrt{\gamma}}|x|}
\end{eqnarray*}
is the fundamental solution to the linear elliptic equation $
\Delta_{x}G-\frac{1}{\gamma}G=0$. Thus \eqref{sta.ellip} admits a
solution if and only if the nonlinear mapping $T$ has a fixed point.
Define
\begin{eqnarray*}
\mathscr{B}_{m,k}(B)=\{\phi\in
W^{m,2}_{k}(\mathbb{R}^{3});\|\phi\|_{W^{m,2}_{k}}\leq
B\|n_{b}-1\|_{W^{m,2}_{k}},\ m\geq2\}
\end{eqnarray*}
for some constant $B$ to be determined later. Next, we prove that if
 $\|n_{b}-1\|_{W^{m,2}_{k}} $ is small enough, there exists a constant $B$ such that
$T:\mathscr{B}_{m,k}(B)\rightarrow \mathscr{B}_{m,k}(B) $ is a
contraction mapping. In fact, for simplicity, let us denote
\begin{eqnarray*}
g(x)=\left(\frac{\gamma-1}{\gamma}x+1\right)^{\frac{1}{\gamma-1}}-\frac{1}{\gamma}x-1.
\end{eqnarray*}
Then it holds that
\begin{eqnarray}\label{T.phi}
T(\phi)(x)=-\int_{\mathbb{R}^{3}}\frac{1}{4\pi|x-y|}e^{-\tfrac{1}{\sqrt{\gamma}}|x-y|}
[g(\phi(y))-(n_{b}(y)-1)]dy.
\end{eqnarray}
Taking derivatives $ \partial_{x}^{\alpha}$ on both sides of $
\eqref{T.phi}$, one has
\begin{eqnarray}\label{T.estimate}
\arraycolsep=1.5pt
\begin{array}[b]{rl}
\partial_{x}^{\alpha}T(\phi)(x)=&\displaystyle-(-1)^{|\alpha|}
\int_{\mathbb{R}^{3}}\frac{1}{4\pi|x-y|}e^{-\tfrac{1}{\sqrt{\gamma}}|x-y|}
[\partial^{\alpha}_{y}
g(\phi(y))-\partial^{\alpha}_{y}(n_{b}(y)-1)]dy\\[5mm]
=&-(-1)^{|\alpha|}G*(\partial^{\alpha}g(\phi)-\partial^{\alpha}(n_{b}-1)).
\end{array}
\end{eqnarray}
Here let's list some properties of the operator $G*$.
\begin{lemma}\label{Pro.OG}
For any $k\geq 0$, it holds that
\begin{eqnarray}\label{pro.Gdecay}
\int_{\mathbb{R}^{3}}\frac{1}{|y|}e^{-\tfrac{1}{\sqrt{\gamma}}|y|}\frac{1}{(1+|x-y|)^{k}}dy\leq
\frac{C_{k}}{(1+|x|)^{k}},
\end{eqnarray}
and for any $f\in W_{k}^{m,2}$,
\begin{eqnarray}\label{pro.G}
\|(1+|x|)^{\frac{k}{2}}(G*f)\|\leq C_{k}^{\frac{1}{2}}\|G\|_{L^{1}}^{\frac{1}{2}}\|(1+|x|)^{\frac{k}{2}}f\|.
\end{eqnarray}
\end{lemma}
\textit{Proof.} \eqref{pro.Gdecay} has been proved in \cite{RY}. We
only prove \eqref{pro.G} by using $\eqref{pro.Gdecay}$.
\begin{eqnarray*}
\arraycolsep=1.5pt
\begin{array}[b]{rl}
\displaystyle\left| \int_{\mathbb{R}^{3}}G(x-y)f(y)dy \right|
\leq & \displaystyle \int_{\mathbb{R}^{3}}\frac{|G(x-y)|^{\frac{1}{2}}}{(1+|y|)^{\frac{k}{2}}}
|G(x-y)|^{\frac{1}{2}}(1+|y|)^{\frac{k}{2}}|f(y)|dy\\[5mm]
\leq &  \displaystyle \left(\int_{\mathbb{R}^{3}}\frac{|G(x-y)|}{(1+|y|)^{k}}dy\right)^{\frac{1}{2}}
\left(\int_{\mathbb{R}^{3}}|G(x-y)|(1+|y|)^{k}|f(y)|^2dy\right)^{\frac{1}{2}}\\[5mm]
\leq &  \displaystyle
\frac{C_{k}^{\frac{1}{2}}}{(1+|x|)^{\frac{k}{2}}}
\left(\int_{\mathbb{R}^{3}}|G(x-y)|(1+|y|)^{k}|f(y)|^2dy\right)^{\frac{1}{2}}.
\end{array}
\end{eqnarray*}
Then
\begin{eqnarray*}
\arraycolsep=1.5pt
\begin{array}[b]{rl}
 & \displaystyle \int_{\mathbb{R}^{3}}(1+|x|)^{k}\left| \int_{\mathbb{R}^{3}}G(x-y)f(y)dy \right|^2dx\\[5mm]
\leq & \displaystyle C_{k} \int_{\mathbb{R}^{3}}\int_{\mathbb{R}^{3}}|G(x-y)|(1+|y|)^{k}|f(y)|^2dydx\\[5mm]
= & \displaystyle  C_{k} \int_{\mathbb{R}^{3}}\int_{\mathbb{R}^{3}}|G(x-y)|(1+|y|)^{k}|f(y)|^2dxdy\\[5mm]
= & \displaystyle  C_{k} \int_{\mathbb{R}^{3}}(1+|y|)^{k}|f(y)|^2 dy \int_{\mathbb{R}^{3}}|G(x-y)|dx\\[5mm]
=& C_{k}\|G\|_{L^{1}}\|(1+|x|)^{\frac{k}{2}}f\|^2.
\end{array}
\end{eqnarray*}
\textbf{Remark:} When $k=0$, $C_{k}=\|G\|_{L^{1}}$, \eqref{pro.G} is  in accordance with Young inequality.

By \eqref{pro.G} and \eqref{T.estimate}, one has
\begin{eqnarray*}
\arraycolsep=1.5pt
\begin{array}[b]{rl}
      & \|(1+|x|)^{\frac{k}{2}}\partial_{x}^{\alpha}T(\phi)(x)\|\\[3mm]
 \leq & C\|(1+|x|)^{\frac{k}{2}}\partial^{\alpha}g(\phi)\|
 +C\|(1+|x|)^{\frac{k}{2}}\partial^{\alpha}(n_{b}-1))\|.
\end{array}
\end{eqnarray*}
By the definition $\eqref{def.norm}$ of the norm $
\|\cdot\|_{W^{m,2}_{k}}$, one has

\begin{eqnarray}\label{T.nb}
\arraycolsep=1.5pt
\begin{array}[b]{rl}
      \|T(\phi)(x)\|_{W^{m,2}_{k}}
   =  & \displaystyle \left(\sum_{|\alpha|\leq m}\|(1+|x|)^{\frac{k}{2}}\partial_{x}^{\alpha}T(\phi)(x)\|^2\right)^{\frac{1}{2}}\\[5mm]
 \leq & \displaystyle C\left(\sum_{|\alpha|\leq m}\|(1+|x|)^{\frac{k}{2}}\partial^{\alpha}g(\phi)\|^2\right)^{\frac{1}{2}}
 +C\left(\sum_{|\alpha|\leq m}\|(1+|x|)^{\frac{k}{2}}\partial^{\alpha}(n_{b}-1)\|^2\right)^{\frac{1}{2}}\\[5mm]
 \leq & \displaystyle C\left(\sum_{|\alpha|\leq m}\|(1+|x|)^{\frac{k}{2}}\partial^{\alpha}g(\phi)\|^2\right)^{\frac{1}{2}}
 +C\|n_{b}-1\|_{W^{m,2}_{k}}.
 \end{array}
\end{eqnarray}
On the other hand, note
\begin{eqnarray*}
 g(\phi)=\int_{0}^{1}\int_{0}^{\theta}g''(\tau\phi)d\tau d\theta
 \phi^{2}\triangleq h(\phi)\phi^{2},
\end{eqnarray*}
where
$g''(x)=\frac{2-\gamma}{\gamma^{2}}\left(\frac{\gamma-1}{\gamma}x+1\right)^{\frac{3-2\gamma}{\gamma-1}}$.

It is straightforward to check that
\begin{eqnarray*}
\|(1+|x|)^{\frac{k}{2}}\partial^{\alpha}(h(\phi)\phi^{2})\|\leq
\sum_{\beta_{1}+\beta_{2}+\beta_{3}=\alpha}C_{\beta_1,\beta_2,\beta_{3}}^{\alpha}
\|(1+|x|)^{\frac{k}{2}}\partial^{\beta_{1}}h(\phi)\partial^{\beta_2}\phi\partial^{\beta_3}\phi\|.
\end{eqnarray*}
 In addition, one has the following claim.

\textbf{Claim}:
\begin{eqnarray}\label{T.gphi}
\|(1+|x|)^{\frac{k}{2}}\partial^{\beta_{1}}h(\phi)\partial^{\beta_2}\phi\partial^{\beta_3}\phi\|\leq
C \|\phi\|^2_{W^{m,2}_{k}}.
\end{eqnarray}

\textit{Proof of claim:} We prove \eqref{T.gphi} by two cases.

\textbf{Case 1.} $\beta_{1}=0$. In this case,
$|\beta_{2}|+|\beta_{3}|\leq m$, thus one can suppose
$|\beta_{2}|\leq [\frac{m}{2}]$ by the symmetry of $\beta_{2}$ and
$\beta_{3}$. This deduces
\begin{eqnarray*}
\arraycolsep=1.5pt
\begin{array}[b]{rl}
     \|(1+|x|)^{\frac{k}{2}}h(\phi)\partial^{\beta_2}\phi\partial^{\beta_3}\phi\|
\leq &
\|h(\phi)\|_{L^{\infty}}\|\partial^{\beta_2}\phi\|_{L^{6}}\|(1+|x|)^{\frac{k}{2}}\partial^{\beta_3}\phi\|_{L^{3}}\\[3mm]
\leq & C \|\nabla
\partial^{\beta_2}\phi\|\|(1+|x|)^{\frac{k}{2}}\partial^{\beta_3}\phi\|_{1}\\[3mm]
\leq & C \|\phi\|^2_{W^{m,2}_{k}}.
\end{array}
\end{eqnarray*}
Here, we have used that $h(\cdot)$ is  a continuous function in the
argument,
\begin{eqnarray*}
\|\phi\|_{L^{\infty}}\leq C\|\nabla \phi\|_{H^{1}}\leq
C\|\phi\|_{W^{m,2}_{k}}\leq C \|n_{b}-1\|_{W^{m,2}_{k}}\ll 1,
\end{eqnarray*}
and $m\geq 2$.

\textbf{Case 2}. $|\beta_{1}|\geq 1$. Notice that
\begin{eqnarray*}
\partial^{\beta_{1}}h(\phi)=\sum_{l=1}^{|\beta_{1}|}h^{(l)}(\phi)
\sum_{\gamma_{1}+\gamma_{2}+\cdots\gamma_{l}=\beta_{1}}
C_{\gamma_{1},\gamma_{2},\cdots\gamma_{l}}\Pi_{i=1}^{l}\partial^{\gamma_{i}}\phi,
\end{eqnarray*}
\eqref{T.gphi} can be similarly obtained because $h^{(m)}(\phi)$ is
also bounded.

Putting $\eqref{T.gphi}$ into \eqref{T.nb}, and using the above
estimates, one has
\begin{eqnarray}\label{T.phi.est}
\|T(\phi)(x)\|_{W^{m,2}_{k}}\leq C B^2
\|n_{b}-1\|_{W^{m,2}_{k}}^2+C\|n_{b}-1\|_{W^{m,2}_{k}}.
\end{eqnarray}

Finally, for any $\phi_{1}=\phi_{1}(x)$ and $\phi_{2}=\phi_{2}(x)$,
it holds that
\begin{eqnarray*}
T(\phi_{1})-T(\phi_{2})=G*(g(\phi_{1})-g(\phi_{2}))
\end{eqnarray*}
with
\begin{eqnarray*}
g(\phi_{1})-g(\phi_{2})=\int_{0}^{1}g'(\theta\phi_{1}+(1-\theta)\phi_{2})d\theta(\phi_{1}-\phi_{2}).
\end{eqnarray*}
Notice that for any $\phi=\phi(x)$,
\begin{eqnarray*}
\arraycolsep=1.5pt
\begin{array}{rcl}
g'(\phi)&=&\displaystyle\frac{1}{\gamma}\left(\frac{\gamma-1}{\gamma}\phi+1\right)^{\frac{2-\gamma}{\gamma-1}}
-\frac{1}{\gamma}\\[3mm]
&=&\displaystyle\int_{0}^{1}\frac{2-\gamma}{\gamma^{2}}
\left(\frac{\gamma-1}{\gamma}\theta\phi+1\right)^{\frac{3-2\gamma}{\gamma-1}}d\theta\phi.
\end{array}
\end{eqnarray*}
Then the same computations as for $\eqref{T.phi.est}$ yield
\begin{eqnarray}\label{T.contract}
\arraycolsep=1.5pt
\begin{array}{rl}
&\|T(\phi_{1})-T(\phi_{2})\|_{W_{k}^{m,2}}\\[3mm]
\leq & C
(\|\phi_{1}\|_{W^{m,2}_{k}}+\|\phi_{2}\|_{W^{m,2}_{k}})\|\phi_{1}-\phi_{2}\|_{W^{m,2}_{k}}.
\end{array}
\end{eqnarray}
Combining $\eqref{T.phi.est}$ with $ \eqref{T.contract}$, the
standard argument implies that $T$ has a unique fixed point $\phi$
in $ \mathscr{B}_{m,k}(B)$ for a proper constant $B$ provided that
$\|n_{b}-1\|_{W^{m,2}_{k}}$ is small enough. This completes Theorem
\ref{sta.existence}.

Let us conclude this section with a remark. The existence of
solutions to the elliptic equation \eqref{sta.ellip} can also be
proved in the framework of \cite{RY} when
$\|n_{b}-1\|_{W^{m,\infty}_{k}}$ is sufficiently small. We consider
the existence  when $\|n_{b}-1\|_{W^{m,2}_{k}}$ is sufficiently
small in order to derive the more general conclusion. In fact, in
the process of dealing with the stability and convergence rates,
only the smallness of $\|n_{b}-1\|_{W^{m,2}_{0}}$ is assumed, and
the space decay at infinity of $n_{b}(x)-1$ is not needed.
\vspace{6mm}

\section{Stability of stationary solution}

\vspace{4mm}
\subsection{Reformulation of the problem} Let $[n,u,E,B]$ be a smooth
solution to the Cauchy problem of the Euler-Maxwell system
(\ref{1.1}) with given initial data (\ref{1.2}) satisfying
(\ref{1.3}). Set

\begin{eqnarray}\label{2.1}
&&\left\{
  \begin{aligned}
   &\sigma(t,x)=\frac{2}{\gamma-1}\left\{\left[n\left(\frac{t}{\sqrt{\gamma}},x\right)\right]
   ^{\frac{\gamma-1}{2}}-1\right\}, \ \ \
   v=\frac{1}{\sqrt{\gamma}}u\left(\frac{t}{\sqrt{\gamma}},x\right),
   \\[5mm]
   &\ \
   \tilde{E}=\frac{1}{\sqrt{\gamma}}E\left(\frac{t}{\sqrt{\gamma}},x\right),\ \
   \ \tilde{B}=\frac{1}{\sqrt{\gamma}}B\left(\frac{t}{\sqrt{\gamma}},x\right).
  \end{aligned}\right.
\end{eqnarray}
Then, $V:=[\sigma,v,\tilde{E},\tilde{B}]$ satisfies

\begin{equation}\label{2.2}
\left\{
  \begin{aligned}
  &\partial_t \sigma+\left(\frac{\gamma-1}{2}\sigma+1\right)\nabla\cdot v+v\cdot \nabla \sigma=0,\\
  &\partial_t v+v \cdot \nabla
  v+\left(\frac{\gamma-1}{2}\sigma+1\right)\nabla \sigma=-\left(\frac{1}{\sqrt{\gamma}}\tilde{E}+v\times \tilde{B}\right)
  -\frac{1}{\sqrt{\gamma}}v,\\
  &\partial_t\tilde{E}-\frac{1}{\sqrt{\gamma}}\nabla\times\tilde{B}
  =\frac{1}{\sqrt{\gamma}}v+\frac{1}{\sqrt{\gamma}}[\Phi(\sigma)+\sigma]v,\\
  &\partial_t \tilde{B}+\frac{1}{\sqrt{\gamma}}\nabla \times \tilde{E}=0,\\
 &\nabla \cdot
  \tilde{E}=-\frac{1}{\sqrt{\gamma}}[\Phi(\sigma)+\sigma]
  +\frac{1}{\sqrt{\gamma}}(n_{b}(x)-1), \ \ \nabla
  \cdot \tilde{B}=0, \ \ \ t>0,\ x\in\mathbb{R}^{3},
\end{aligned}\right.
\end{equation}
with initial data
\begin{eqnarray}\label{2.3}
V|_{t=0}=V_{0}:=[\sigma_{0},v_{0},\tilde{E}_{0},\tilde{B}_{0}],\ \
x\in\mathbb{R}^{3}.
\end{eqnarray}
Here, $\Phi(\cdot)$  is defined by
\begin{eqnarray}\label{def.phi}
\Phi(\sigma)=\left(\frac{\gamma-1}{2}\sigma+1\right)^{\frac{2}{\gamma-1}}-\sigma-1,
\end{eqnarray}
and $V_{0}=[\sigma_{0},v_{0},\tilde{E}_{0},\tilde{B}_{0}]$ is given
from $[n_{0},u_{0},E_{0},B_0]$ according to the transform
(\ref{2.1}), and hence $V_{0}$ satisfies
\begin{eqnarray}\label{2.4}
   \nabla\cdot\tilde{E}_0=-\frac{1}{\sqrt{\gamma}}[\Phi(\sigma_{0})+\sigma_{0}]
   +\frac{1}{\sqrt{\gamma}}(n_{b}(x)-1),\ \ \ \
   \nabla \cdot \tilde{B}_0=0,\ \ \ x\in\mathbb{R}^{3}.
\end{eqnarray}

On the other hand, set
\begin{eqnarray}\label{sta.tran}
 \sigma_{st}(x)=\frac{2}{\gamma-1}\left\{n_{st}(x)^{\frac{\gamma-1}{2}}-1\right\},
 \ \ \ \
  \tilde{E}_{st}=\frac{1}{\sqrt{\gamma}}E_{st}(x).
\end{eqnarray}
Then, $[\sigma_{st},\tilde{E}_{st}]$ satisfies

\begin{eqnarray}\label{sta.eq}
\left\{\begin{aligned}
& \left(\frac{\gamma-1}{2}\sigma_{st}+1\right)\nabla\sigma_{st}=-\frac{1}{\sqrt{\gamma}}\tilde{E}_{st},\\
&\frac{1}{\sqrt{\gamma}}\nabla\times \tilde{E}_{st}=0,\\
 &\nabla \cdot
 \tilde{E}_{st}=\frac{1}{\sqrt{\gamma}}(n_{b}(x)-1)-\frac{1}{\sqrt{\gamma}}(\Phi(\sigma_{st})+\sigma_{st}).
\end{aligned}\right.
\end{eqnarray}
Based on the existence result proved in Section 2, we will study the
stability of the stationary state
$[\sigma_{st},0,\tilde{E}_{st},0]$. Set the perturbations
$[\bar{\sigma},\bar{v},\bar{E},\bar{B}]$ by
\begin{eqnarray*}
\bar{\sigma}=\sigma-\sigma_{st},\ \ \bar{v}=v,\ \
\bar{E}=\tilde{E}-\tilde{E}_{st},\ \  \bar{B}=\tilde{B}.
\end{eqnarray*}
Combining \eqref{2.2} with \eqref{sta.eq}, then
$\bar{V}:=[\bar{\sigma},\bar{v},\bar{E},\bar{B}]$ satisfies

\begin{equation}\label{sta.equ}
\left\{
 \begin{aligned}
  &\partial_t \bar{\sigma}+(\frac{\gamma-1}{2}\bar{\sigma}+1)\nabla\cdot \bar{v}+\bar{v}\cdot \nabla \bar{\sigma}
  +\bar{v}\cdot \nabla \sigma_{st}+\frac{\gamma-1}{2}\sigma_{st}\nabla\cdot \bar{v}=0,\\
  &\partial_t \bar{v}+\bar{v} \cdot \nabla
  \bar{v}+(\frac{\gamma-1}{2}\bar{\sigma}+1)\nabla \bar{\sigma}
  +\frac{\gamma-1}{2}\bar{\sigma}\nabla \sigma_{st}+\frac{\gamma-1}{2}\sigma_{st}\nabla\bar{\sigma}=
  -(\frac{1}{\sqrt{\gamma}}\bar{E}+\bar{v}\times \bar{B})
  -\frac{1}{\sqrt{\gamma}}\bar{v},\\
  &\partial_t\bar{E}-\frac{1}{\sqrt{\gamma}}\nabla\times \bar{B}
  =\frac{1}{\sqrt{\gamma}}\bar{v}+
  \frac{1}{\sqrt{\gamma}}[\Phi(\bar{\sigma}+\sigma_{st})+\bar{\sigma}+\sigma_{st}]\bar{v},\\
  &\partial_t \bar{B}+\frac{1}{\sqrt{\gamma}}\nabla \times \bar{E}=0,\\
 &\nabla \cdot
  \bar{E}=-\frac{1}{\sqrt{\gamma}}[\Phi(\bar{\sigma}+\sigma_{st})-\Phi(\sigma_{st})]
  -\frac{1}{\sqrt{\gamma}}\bar{\sigma}, \ \ \nabla
  \cdot \bar{B}=0,  \ \ t>0,\ x\in\mathbb{R}^{3},
\end{aligned}\right.
\end{equation}
with initial data
\begin{eqnarray}\label{sta.equi}
\bar{V}|_{t=0}=\bar{V}_{0}:=[\sigma_{0}-\sigma_{st},v_{0},\tilde{E}_{0}-\tilde{E}_{st},\tilde{B}_{0}],\
\ x\in\mathbb{R}^{3}.
\end{eqnarray}
Here, $\Phi(\cdot)$  is defined by \eqref{def.phi},
 and  $\bar{V}_{0}$ satisfies
\begin{eqnarray}\label{sta.equC}
   \nabla \cdot
  \bar{E}_{0}=-\frac{1}{\sqrt{\gamma}}[\Phi(\bar{\sigma}_{0}+\sigma_{st})-\Phi(\sigma_{st})]
  -\frac{1}{\sqrt{\gamma}}\bar{\sigma}_{0}, \ \ \nabla
  \cdot \bar{B}_{0}=0, \ \ t>0,\ x\in\mathbb{R}^{3}.
\end{eqnarray}

In what follows, we suppose the integer $N \geq 3$. Besides, for
$\bar{V}=[\bar{\sigma},\bar{v},\bar{E},\bar{B}]$, we define the full
instant energy functional $\mathcal {E}_{N}(\bar{V}(t))$, the
high-order instant energy functional $\mathcal
{E}_{N}^{h}(\bar{V}(t))$, and the dissipation rates $\mathcal
{D}_{N}(\bar{V}(t))$, $\mathcal {D}_{N}^{h}(\bar{V}(t))$ by

\begin{equation}\label{de.E}
\arraycolsep=1.5pt
\begin{array}{rl}
\mathcal{E}_{N}(\bar{V}(t))=&\displaystyle\sum_{|\alpha|\leq
N}\int_{\mathbb{R}^3}(1+\sigma_{st}+\Phi(\sigma_{st}))
(|\partial^{\alpha}\bar{\sigma}|^2+|\partial^{\alpha}\bar{v}|^2)dx+\|[\bar{E},\bar{B}]\|_{N}^{2}\\[5mm]
&\displaystyle+\kappa_{1}\sum_{|\alpha|\leq N-1} \langle
\partial^{\alpha}\bar{v},\nabla\partial^{\alpha}\bar{\sigma}\rangle+\kappa_{2}\sum_{|\alpha|\leq N-1}\langle
\partial^{\alpha}\bar{v},\partial^{\alpha}\bar{E}\rangle\\[5mm]
&\displaystyle-\kappa_{3}\sum_{|\alpha|\leq N-2}\langle \nabla
\times
\partial^{\alpha}\bar{E},\partial^{\alpha}\bar{B}\rangle,
\end{array}
\end{equation}
and
\begin{equation}\label{de.Eh}
\begin{aligned}
\mathcal{E}_{N}^{h}(\bar{V}(t))&=\sum_{1\leq|\alpha|\leq
N}\int_{\mathbb{R}^3}(1+\sigma_{st}+\Phi(\sigma_{st}))
(|\partial^{\alpha}\bar{\sigma}|^2+|\partial^{\alpha}\bar{v}|^2)dx+\|\nabla[\bar{E},\bar{B}]\|_{N-1}^{2}\\
&+\kappa_{1}\sum_{1\leq|\alpha|\leq N-1}\langle
\partial^{\alpha}\bar{v},\nabla\partial^{\alpha}\bar{\sigma}\rangle+\kappa_{2}\sum_{1\leq|\alpha|\leq N-1}\langle
\partial^{\alpha}\bar{v},\partial^{\alpha}\bar{E}\rangle\\[3mm]
&-\kappa_{3}\sum_{1\leq |\alpha|\leq N-2}\langle \nabla
   \times\partial^{\alpha}\bar{E},\partial^{\alpha}\bar{B}\rangle,
\end{aligned}
\end{equation}
respectively, where $0<\kappa_{3}\ll\kappa_{2}\ll\kappa_{1}\ll 1$
are constants to be properly chosen in the later proof. Notice that
since all constants $\kappa_i$ $(i=1,2,3)$ are small enough, one has
\begin{equation*}
    \mathcal {E}_{N}(\bar{V}(t))\sim
\|[\bar{\sigma},\bar{v},\bar{E},\bar{B}] \|_{N}^{2},\quad \mathcal
{E}_{N}^{h}(\bar{V}(t))\sim \|\nabla
[\bar{\sigma},\bar{v},\bar{E},\bar{B}]  \|_{N-1}^{2}.
\end{equation*}
We further define the dissipation rates $\mathcal
{D}_{N}(\bar{V}(t))$, $\mathcal {D}_{N}^{h}(\bar{V}(t))$ by

\begin{eqnarray}\label{de.D}
\arraycolsep=1.5pt
\begin{array}{rl}
\mathcal {D}_{N}(\bar{V}(t))=\displaystyle \sum_{|\alpha|\leq
N}\int_{\mathbb{R}^3}(1+\sigma_{st}&+\Phi(\sigma_{st}))|\partial^{\alpha}\bar{v}|^{2}dx\\[3mm]
&+\|\bar{\sigma}\|_{N}^{2}+\|\nabla[\bar{E},\bar{B}]\|_{N-2}^{2}+\|\bar{E}\|^{2},
\end{array}
\end{eqnarray}
and
\begin{eqnarray}\label{de.Dh}
\arraycolsep=1.5pt
\begin{array}{rl}
\mathcal {D}_{N}^{h}(\bar{V}(t))=\displaystyle
\sum_{1\leq|\alpha|\leq
N}\int_{\mathbb{R}^3}(1+\sigma_{st}&+\Phi(\sigma_{st}))|\partial^{\alpha}\bar{v}|^{2}dx\\[3mm]
&+\|\nabla\bar{\sigma}\|_{N-1}^{2}+\|\nabla^2[\bar{E},\bar{B}]\|_{N-3}^{2}+\|\nabla\bar{E}\|^{2}.
\end{array}
\end{eqnarray}
Then, concerning the reformulated Cauchy problem
$\eqref{sta.equ}$-$\eqref{sta.equi}$, one has the following global
existence result.

\begin{proposition}\label{pro.2.1}
Suppose that $\|n_{b}-1\|_{W_{0}^{N+1,2}}$ is small enough and
$\eqref{sta.equC}$ holds for given initial data
$\bar{V}_{0}=[\sigma_{0}-\sigma_{st},v_{0},\tilde{E}_0-\tilde{E}_{st},\tilde{B}_{0}]$.
Then, there are $\mathcal {E}_{N}(\cdot) $ and $\mathcal
{D}_{N}(\cdot)$ in the form $\eqref{de.E} $ and $\eqref{de.D}$ such
that the following holds true:

If $\mathcal {E}_{N}(\bar{V}_{0})>0$ is small enough,  the Cauchy
problem $\eqref{sta.equ}$-$\eqref{sta.equi}$ admits a unique global
nonzero solution
$\bar{V}=[\sigma-\sigma_{st},v,\tilde{E}-\tilde{E}_{st},\tilde{B}] $
satisfying
\begin{eqnarray}\label{V.satisfy}
\bar{V} \in C([0,\infty);H^{N}(\mathbb{R}^{3}))\cap {\rm
Lip}([0,\infty);H^{N-1}(\mathbb{R}^{3})),
\end{eqnarray}
and
\begin{eqnarray}\label{pro.2.1j}
\mathcal {E}_{N}(\bar{V}(t))+\lambda\int_{0}^{t}\mathcal
{D}_{N}(\bar{V}(s))ds\leq \mathcal {E}_{N}(\bar{V}_{0})
\end{eqnarray}
for any $t\geq 0$.
\end{proposition}

Moreover, solutions obtained in Proposition $ \ref{pro.2.1}$ indeed
decay in time with some rates under some extra regularity and
integrability conditions on initial data. For that, given
$\bar{V}_{0}=[\sigma_{0}-\sigma_{st},v_{0},\tilde{E}_0-\tilde{E}_{st},\tilde{B}_{0}]$,
set $\epsilon_{m}(\bar{V}_0)$ as
\begin{eqnarray}\label{def.epsi}
\epsilon_{m}(\bar{V}_0)=\|\bar{V}_{0}\|_{m}+\|[v_{0},\tilde{E}_0-\tilde{E}_{st},\tilde{B}_{0}]\|_{L^{1}},
\end{eqnarray}
for the  integer $m \geq 6$. Then one has the following proposition.
\begin{proposition}\label{pro.2.2}
Suppose that $\|n_{b}-1\|_{W_{0}^{N+4,2}}$ is small enough and
$\eqref{sta.equC}$ holds for given initial data
$\bar{V}_{0}=[\sigma_{0}-\sigma_{st},v_{0},\tilde{E}_0-\tilde{E}_{st},\tilde{B}_{0}]$.
If $\epsilon_{N+3}(\bar{V}_{0})>0$ is small enough, then the
solution
$\bar{V}=[\sigma-\sigma_{st},v,\tilde{E}-\tilde{E}_{st},\tilde{B}] $
satisfies
\begin{eqnarray}\label{V.decay}
\|\bar{V}(t)\|_{N} \leq C
\epsilon_{N+3}(\bar{V}_{0})(1+t)^{-\frac{3}{4}},
\end{eqnarray}
and
\begin{eqnarray}\label{nablaV.decay}
\|\nabla \bar{V}(t)\|_{N-1} \leq C
\epsilon_{N+3}(\bar{V}_{0})(1+t)^{-\frac{5}{4}}
\end{eqnarray}
for any $t\geq 0$.
\end{proposition}

\subsection{a priori estimates} In this subsection, we prove that the stationary
solution obtained in Section \ref{sec2} is stable under small
initial perturbation.
 We begin to use the refined energy
method to obtain some uniform-in-time {\it a priori} estimates for
smooth solutions to the Cauchy problem
(\ref{sta.equ})-(\ref{sta.equi}). To the end, let us denote
\begin{equation}\label{def.delta}
\delta=\|\sigma_{st}\|_{W_{0}^{N+1,2}}=\left(\sum_{|\alpha|\leq
N+1}\int_{\mathbb{R}^3}|\partial^{\alpha}_{x}\sigma_{st}|^2dx\right)^{\frac{1}{2}}
\end{equation}
for simplicity of presentation. A careful look at the proof of
Theorem \ref{sta.existence} shows that
\begin{eqnarray*}
\sigma_{st}&=&\frac{2}{\gamma-1}\left\{n_{st}^{\frac{\gamma-1}{2}}-1\right\}\\
                &=&\frac{2}{\gamma-1}\left\{\left(\frac{\gamma-1}{\gamma}Q_{st}+1\right)^{\frac{1}{2}}-1\right\}\\
                &=&\frac{2}{\gamma}\dfrac{Q_{st}}{\left(\frac{\gamma-1}{\gamma}Q_{st}+1\right)^{\frac{1}{2}}+1}\sim
                Q_{st}.
\end{eqnarray*}
It follows that $\delta \leq C\|Q_{st}\|_{W_{0}^{N+1,2}}\leq
C\|n_{b}-1\|_{W_{0}^{N+1,2}}$ is small enough. Notice that
(\ref{sta.equ}) is a quasi-linear symmetric hyperbolic system. The
main goal of this subsection is to prove

\begin{theorem}\label{estimate}(\textrm{a priori estimates}).
Let $0<T\leq \infty$ be given. Suppose
$\bar{V}:=[\bar{\sigma},\bar{v},\bar{E},\bar{B}]\in
C([0,T);H^{N}(\mathbb{R}^{3}))$ is smooth for $T>0$ with
\begin{eqnarray}\label{3.1}
\sup_{0\leq t<T}\|\bar{V}(t)\|_{N}\leq 1,
\end{eqnarray}
and assume that $\bar{V}$ solves the system (\ref{sta.equ}) for
$t\in(0,T)$. Then, there are $\mathcal {E}_{N}(\cdot) $ and
$\mathcal {D}_{N}(\cdot)$ in the form $\eqref{de.E} $ and
$\eqref{de.D}$ such that
\begin{eqnarray}\label{3.2}
&& \frac{d}{dt}\mathcal {E}_{N}(\bar{V}(t))+\lambda\mathcal
{D}_{N}(\bar{V}(t))\leq
 C[\mathcal {E}_{N}(\bar{V}(t))^{\frac{1}{2}}+\mathcal {E}_{N}(\bar{V}(t))+\delta]\mathcal {D}_{N}(\bar{V}(t))
\end{eqnarray}
for any $0\leq t<T$.
\end{theorem}

\begin{proof}
The proof is divided into five steps.

\medskip
 \textbf{ Step 1.} It holds that
\begin{equation}\label{3.3}
\begin{aligned}
&\frac{1}{2}\frac{d}{dt}\left(\sum_{|\alpha|\leq
N}\int_{\mathbb{R}^3}(1+\sigma_{st}+\Phi(\sigma_{st}))
(|\partial^{\alpha}\bar{\sigma}|^2+|\partial^{\alpha}\bar{v}|^2)dx+\|[\bar{E},\bar{B}]\|_{N}^{2}\right)\\
&+\frac{1}{\sqrt{\gamma}}\sum_{|\alpha|\leq
N}\int_{\mathbb{R}^3}(1+\sigma_{st}+\Phi(\sigma_{st}))|\partial^{\alpha}\bar{v}|^{2}dx\\
\leq &
C(\|\bar{V}\|_{N}+\delta)(\|[\bar{\sigma},\bar{v}]\|^{2}+\|\nabla[\bar{\sigma},\bar{v}]\|_{N-1}^{2}
 +\|\nabla \bar{E}\|_{N-2}^2).
 \end{aligned}
\end{equation}
In fact, applying $\partial^{\alpha}$ to the first two equations of
(\ref{sta.equ}) for $|\alpha|\leq N$ and multiplying them by
$(1+\sigma_{st}+\Phi(\sigma_{st}))\partial^{\alpha}\bar{\sigma}$ and
$(1+\sigma_{st}+\Phi(\sigma_{st}))\partial^{\alpha}\bar{v}$
respectively, taking integrations in $x$ and then using integration
by parts give

\begin{eqnarray}\label{3.4}
  && \begin{aligned}
   &\frac{1}{2}\frac{d}{dt}\int_{\mathbb{R}^3}(1+\sigma_{st}+\Phi(\sigma_{st}))
   (|\partial^{\alpha}\bar{\sigma}|^2+|\partial^{\alpha}\bar{v}|^2)dx+\frac{1}{\sqrt{\gamma}}
   \langle \partial^{\alpha}\bar{E},(1+\sigma_{st}+\Phi(\sigma_{st}))\partial^{\alpha}\bar{v}\rangle\\
   &+\frac{1}{\sqrt{\gamma}}\int_{\mathbb{R}^3}(1+\sigma_{st}+\Phi(\sigma_{st}))|\partial^{\alpha}\bar{v}|^{2}dx
   =-\sum_{\beta<\alpha}C^{\alpha}_{\beta}I_{\alpha,\beta}(t)+I_{1}(t).
 \end{aligned}
\end{eqnarray}
Here,
$I_{\alpha,\beta}(t)=I_{\alpha,\beta}^{(\sigma)}(t)+I_{\alpha,\beta}^{(v)}(t)$
with
\begin{eqnarray*}
\arraycolsep=1.5pt
\begin{array}{rl}
 \displaystyle
I_{\alpha,\beta}^{(\sigma)}(t)=& \displaystyle\langle
\partial^{\alpha-\beta}\bar{v} \cdot \nabla \partial^{\beta}\bar{\sigma},
(1+\sigma_{st}+\Phi(\sigma_{st}))\partial^{\alpha}\bar{\sigma}\rangle\\[3mm]
&\displaystyle+\frac{\gamma-1}{2}\langle\partial^{\alpha-\beta}\bar{\sigma}
\partial^{\beta}\nabla\cdot\bar{v} ,(1+\sigma_{st}+\Phi(\sigma_{st}))\partial^{\alpha}\bar{\sigma}\rangle\\[3mm]
 & \displaystyle+\frac{\gamma-1}{2}\langle
\partial^{\alpha-\beta}\sigma_{st}  \partial^{\beta}\nabla\cdot\bar{v} ,
(1+\sigma_{st}+\Phi(\sigma_{st}))\partial^{\alpha}\bar{\sigma}
\rangle\\[3mm]
& \displaystyle +\langle
\partial^{\alpha-\beta}\bar{v}\cdot \partial^{\beta}\nabla \sigma_{st} ,(1+\sigma_{st}+\Phi(\sigma_{st}))\partial^{\alpha}\bar{\sigma}\rangle,
 \end{array}
\end{eqnarray*}
\begin{eqnarray*}
\arraycolsep=1.5pt
\begin{array}{rl}
 \displaystyle
I_{\alpha,\beta}^{(v)}(t)=& \displaystyle\langle
\partial^{\alpha-\beta}\bar{v} \cdot \nabla \partial^{\beta}\bar{v} ,
(1+\sigma_{st}+\Phi(\sigma_{st}))\partial^{\alpha}\bar{v}\rangle\\[3mm]
&\displaystyle+\frac{\gamma-1}{2}\langle\partial^{\alpha-\beta}\bar{\sigma}
\nabla\partial^{\beta}\bar{\sigma},(1+\sigma_{st}+\Phi(\sigma_{st}))\partial^{\alpha}\bar{v}\rangle\\[3mm]
&\displaystyle+\frac{\gamma-1}{2}\langle
\partial^{\alpha-\beta}\sigma_{st} \nabla \partial^{\beta}\bar{\sigma}
,(1+\sigma_{st}+\Phi(\sigma_{st}))\partial^{\alpha}\bar{v}\rangle\\[3mm]
&\displaystyle+\langle
\partial^{\alpha-\beta}\bar{v}\times  \partial^{\beta}\bar{B}
,(1+\sigma_{st}+\Phi(\sigma_{st}))\partial^{\alpha}\bar{v}\rangle\\[3mm]
&\displaystyle+\frac{\gamma-1}{2}\langle
\partial^{\alpha-\beta}\bar{\sigma} \nabla \partial^{\beta}\sigma_{st} ,(1+\sigma_{st}+\Phi(\sigma_{st}))
\partial^{\alpha}\bar{v} \rangle
\end{array}
\end{eqnarray*}
and
\begin{eqnarray*}
\arraycolsep=1.5pt
\begin{array}{rl}
 I_{1}(t)=&
 \displaystyle\frac{1}{2}\langle \nabla \cdot \bar{v},
 (1+\sigma_{st}+\Phi(\sigma_{st}))(|\partial^{\alpha}\bar{\sigma}|^{2}+|\partial^{\alpha}\bar{v}|^{2})
 \rangle\\[3mm]
 & \displaystyle+ \frac{\gamma-1}{2}\langle \nabla
\bar{\sigma}\cdot\partial^{\alpha}\bar{v}
,(1+\sigma_{st}+\Phi(\sigma_{st}))\partial^{\alpha}\bar{\sigma}
\rangle-\langle \bar{v}\times
\partial^{\alpha}\bar{B},(1+\sigma_{st}+\Phi(\sigma_{st}))\partial^{\alpha}\bar{v}\rangle\\[3mm]
& \displaystyle+\frac{\gamma-1}{2}\langle \nabla\sigma_{st}
\partial^{\alpha}\bar{v},(1+\sigma_{st}+\Phi(\sigma_{st}))\partial^{\alpha}\bar{\sigma}
\rangle-\frac{\gamma-1}{2}\langle \bar{\sigma}\partial^{\alpha}
\nabla\sigma_{st},(1+\sigma_{st}+\Phi(\sigma_{st}))\partial^{\alpha}\bar{v}
\rangle\\[3mm]
&- \displaystyle \langle \bar{v}\cdot\partial^{\alpha}
\nabla\sigma_{st},(1+\sigma_{st}+\Phi(\sigma_{st}))\partial^{\alpha}\bar{\sigma}
\rangle\\[3mm]
&\displaystyle+\left\langle
\left(\frac{\gamma-1}{2}\bar{\sigma}+1\right)\partial^{\alpha}\bar{v},
\nabla(1+\sigma_{st}+\Phi(\sigma_{st}))\partial^{\alpha}\bar{\sigma}
\right\rangle\\[3mm]
& \displaystyle +\frac{\gamma-1}{2}\langle \sigma_{st}
\partial^{\alpha}\bar{v},\nabla(1+\sigma_{st}+\Phi(\sigma_{st}))\partial^{\alpha}\bar{\sigma}
\rangle\\[3mm]
&\displaystyle+\frac{1}{2}\langle\bar{v},\nabla(1+\sigma_{st}+\Phi(\sigma_{st}))
(|\partial^{\alpha}\bar{\sigma}|^2+|\partial^{\alpha}\bar{v}|^{2})\rangle\triangleq\sum_{j=1}^{9}I_{1,j}(t).
\end{array}
\end{eqnarray*}
When $|\alpha|=0$, it suffices to estimate $I_{1}(t)$ by
\begin{eqnarray*}
\begin{aligned}
 I_{1}(t)\leq &C \|\nabla \cdot \bar{v}\|(\|\bar{v}\|_{L^{6}}\|\bar{v}\|_{L^{3}}
 +\|\bar{\sigma}\|_{L^{6}}\|\bar{\sigma}\|_{L^{3}})
 +C \|\nabla \bar{\sigma}\|\|\bar{v}\|_{L^{6}}\|\bar{\sigma}\|_{L^{3}}
 +C \|\bar{B}\|_{L^{\infty}}\|\bar{v}\|^{2}\\
  &+C \|\nabla\sigma_{st}\|\left\|\bar{\sigma}\right\|_{L^{6}}
 \|\bar{v}\|_{L^3}+C\|\sigma_{st}\|_{L^{\infty}}\|\nabla\sigma_{st}\|\left\|\bar{\sigma}\right\|_{L^{6}}
 \|\bar{v}\|_{L^3}\\
  &+\|\bar{v}\|_{L^{\infty}}
  \|\nabla\sigma_{st}\|(\|\bar{\sigma}\|_{L^{6}}\|\bar{\sigma}\|_{L^{3}}
  +\|\bar{v}\|_{L^{6}}\|\bar{v}\|_{L^{3}})
  \\
  \leq & C (\|[\bar{\sigma},\bar{v}]\|_{H^{1}}+\delta+\delta\|\nabla\bar{v}\|_{H^1})(\|\nabla
  [\bar{\sigma},\bar{v}]\|^{2}+\|[\bar{\sigma},\bar{v}]\|^{2})+ C \|\nabla
  \bar{B}\|_{H^{1}}\|\bar{v}\|^{2},
 \end{aligned}
\end{eqnarray*}
which is further bounded by the r.h.s. term of (\ref{3.3}). When
$|\alpha|\geq 1$, for $I_{1}(t)$, the similarity of $I_{1,1}(t)$ and
$I_{1,2}(t)$ shows that we can estimate them together as follows
\begin{eqnarray*}
\begin{aligned}
I_{1,1}(t)+I_{1,2}(t)\leq & C
\|\nabla\cdot\bar{v}\|_{L^{\infty}}\|(1+\sigma_{st}+\Phi(\sigma_{st}))\|_{L^{\infty}}
\|\nabla[\bar{\sigma},\bar{v}]\|_{N-1}^2\\
&+C\|\nabla\bar{\sigma}\|_{L^{\infty}}\|(1+\sigma_{st}+\Phi(\sigma_{st}))\|_{L^{\infty}}
\|\nabla[\bar{\sigma},\bar{v}]\|_{N-1}^2\\
\leq & C
\|[\bar{\sigma},\bar{v}]\|_{N}\|\nabla[\bar{\sigma},\bar{v}]\|_{N-1}^2.
 \end{aligned}
\end{eqnarray*}
For $I_{1,3}(t)$, $I_{1,5}(t)$ and $I_{1,6}(t)$, there are no
derivative of $\bar{\sigma}$ or $\bar{v}$, then we use $L^{\infty}$
of $\bar{v}$ or $\bar{\sigma}$,
\begin{eqnarray*}
\begin{aligned}
I_{1,3}(t)+I_{1,5}(t)+I_{1,6}(t)\leq & C
\|\bar{v}\|_{L^{\infty}}\|\partial^{\alpha}
\bar{B}\|\|(1+\sigma_{st}+\Phi(\sigma_{st}))\|_{L^{\infty}}
\|\nabla\bar{v}\|_{N-1}\\
&+C \|\bar{\sigma}\|_{L^{\infty}}\|\partial^{\alpha}
\nabla\sigma_{st}\|\|(1+\sigma_{st}+\Phi(\sigma_{st}))\|_{L^{\infty}}
\|\nabla\bar{v}\|_{N-1}\\
&+C \|\bar{v}\|_{L^{\infty}}\|\partial^{\alpha}
\nabla\sigma_{st}\|\|(1+\sigma_{st}+\Phi(\sigma_{st}))\|_{L^{\infty}}
\|\nabla\bar{\sigma}\|_{N-1}\\
\leq &
C(\delta+\|\bar{B}\|_{N})\|\nabla[\bar{\sigma},\bar{v}]\|_{N-1}^2.
 \end{aligned}
\end{eqnarray*}
For other terms of $I_{1}(t)$, both $\bar{\sigma}$ and $\bar{v}$
contain the derivative, one can use the $L^{2}$ of these terms and
$L^{\infty}$ of others. Combining the above two estimates, one has
\begin{eqnarray*}
I_{1}(t)\leq C (\|[\bar{\sigma},\bar{v},\bar{B}]\|_{N}
  +\delta+\delta\|\nabla \bar{v}\|_{H^1})\|\nabla
  [\bar{\sigma},\bar{v}]\|_{N-1}^{2},
\end{eqnarray*}
which is bounded by the r.h.s. term of (\ref{3.3}). On the other
hand, since each term in $I_{\alpha,\beta}(t)$ is the integration of
the four-terms product in which there is at least one term
containing the derivative, one has
\begin{eqnarray*}
I_{\alpha,\beta}(t)\leq C (\|[\bar{\sigma},\bar{v},\bar{B}]\|_{N}
  +\delta+\delta\|\nabla \bar{v}\|_{H^1})\|\nabla
  [\bar{\sigma},\bar{v}]\|_{N-1}^{2},
\end{eqnarray*}
which is also bounded by the r.h.s. term of (\ref{3.3}).

From (\ref{sta.equ}), energy estimates on $\partial^{\alpha}\bar{E}$
and $\partial^{\alpha}\bar{B}$ with $|\alpha| \leq N$ give
\begin{eqnarray}\label{3.5}
 && \begin{aligned}
   &\frac{1}{2}\frac{d}{dt}\|\partial^{\alpha}[\bar{E},\bar{B}]\|^{2}
-\frac{1}{\sqrt{\gamma}}\langle
(1+\sigma_{st}+\Phi(\sigma_{st}))\partial
   ^{\alpha}\bar{v},\partial^{\alpha}\bar{E}\rangle\\
   =&\frac{1}{\sqrt{\gamma}}\langle \partial
   ^{\alpha}[(\Phi(\bar{\sigma}+\sigma_{st})-\Phi(\sigma_{st}))\bar{v}],\partial^{\alpha}\bar{E}\rangle+
   \frac{1}{\sqrt{\gamma}}\langle \partial
   ^{\alpha}[\bar{\sigma}\bar{v}],\partial^{\alpha}\bar{E}\rangle\\
   &+\frac{1}{\sqrt{\gamma}}\sum_{\beta<\alpha}C_{\beta}^{\alpha}\langle \partial^{\alpha-\beta}(1+\sigma_{st}+\Phi(\sigma_{st}))
   \partial^{\beta}\bar{v},\partial^{\alpha}\bar{E}\rangle\\
   =&I_{2,1}(t)+I_{2,2}(t)+\sum_{\beta<\alpha}C_{\beta}^{\alpha}I_{2,\beta}(t).
 \end{aligned}
\end{eqnarray}
In a similar way as before,  when $|\alpha|=0$, it suffices to
estimate $I_{2,1}(t)+ I_{2,2}(t)$ by
\begin{eqnarray*}
  I_{2,1}(t)+ I_{2,2}(t)\leq C \|\nabla
  \bar{\sigma}\|\cdot\|\bar{v}\|_{1}\|\bar{E}\|.
\end{eqnarray*}
When $|\alpha|>0$, $I_{2,1}(t)$ and $I_{2,2}(t)$ can be estimated in
a similar way as in \cite{Duan},
\begin{eqnarray*}
  I_{2,1}(t)+ I_{2,2}(t)\leq C \|\nabla \bar{\sigma}\|_{N-1}\|\nabla \bar{v}\|_{N-1}\|\bar{E}\|_{N}.
\end{eqnarray*}
When $|\alpha|>0$, for each $\beta$ with $\beta<\alpha$,
$I_{2,\beta}$ is estimated by three cases.

\textsl{Case 1.} $|\alpha|=N$. In this case, integration by parts
shows that
\begin{eqnarray*}
 && \begin{aligned}
  I_{2,\beta}(t) \leq & C \delta \|\nabla \bar{v}\|_{N-1}\|\nabla\bar{E}\|_{N-2}\\
  \leq &  C \delta \|\nabla \bar{v}\|_{N-1}^2+C \delta \|\nabla \bar{E}\|_{N-2}^2.
\end{aligned}
\end{eqnarray*}

\textsl{Case 2.} $|\alpha|<N $ and $|\beta|\geq 1$ which imply
$|\alpha-\beta|\leq N-2$. It holds that
\begin{eqnarray*}
 && \begin{aligned}
  I_{2,\beta}(t) \leq &  C\|\partial^{\alpha-\beta}(1+\sigma_{st}+\Phi(\sigma_{st}))\|_{L^{\infty}}
\|\partial^{\beta}\bar{v}\|\|\partial^{\alpha}\bar{E}\|\\
\leq & C\|\nabla\partial^{\alpha-\beta}(1+\sigma_{st}+\Phi(\sigma_{st}))\|_{H^{1}}\|\nabla \bar{v}\|_{N-1}\|\nabla\bar{E}\|_{N-2}\\
  \leq &  C \delta \|\nabla \bar{v}\|_{N-1}^2+C \delta \|\nabla
  \bar{E}\|_{N-2}^2.
\end{aligned}
\end{eqnarray*}

\textsl{Case 3.} $|\alpha|<N $ and $|\beta|=0$. In this case, there
is no derivative of $\bar{v}$, one can use $L^{\infty}$ of $\bar{v}$
to estimate $I_{2,\beta}(t)$,
\begin{eqnarray*}
 && \begin{aligned}
  I_{2,\beta}(t) \leq &  C\|\partial^{\alpha-\beta}(1+\sigma_{st}+\Phi(\sigma_{st}))\|
\|\bar{v}\|_{L^{\infty}}\|\partial^{\alpha}\bar{E}\|\\
  \leq &  C \delta \|\nabla \bar{v}\|_{N-1}^2+C \delta \|\nabla
  \bar{E}\|_{N-2}^2,
\end{aligned}
\end{eqnarray*}
which is bounded by the r.h.s. term of (\ref{3.3}). Then (\ref{3.3})
follows by taking summation of (\ref{3.4}) and (\ref{3.5}) over
$|\alpha| \leq N$. Then the time evolution of the full instant
energy $\|V(t)\|_{N}^{2}$ has been obtained but its dissipation rate
only contains the contribution from the explicit relaxation variable
$\bar{v}$. In a parallel way as \cite{Duan}, by introducing some
interactive functionals, the dissipation from contributions of the
rest components $\bar{\sigma}$, $\bar{E}$, and $\bar{B}$ can be
recovered in turn.

\medskip

\textbf{Step 2.} It holds that
\begin{eqnarray}\label{step2}
  &&\begin{aligned}
  &\frac{d}{dt}\mathcal {E}_{N,1}^{int}(\bar{V})+\lambda\|\bar{\sigma}\|^{2}_{N} \\
  \leq & C\|\nabla\bar{v}\|_{N-1}^{2}+C(\|[\bar{\sigma}, \bar{v},\bar{B}]\|_{N}^{2}+\delta)
  \|\nabla[\bar{\sigma},\bar{v}]\|_{N-1}^{2},
  \end{aligned}
\end{eqnarray}
where $\mathcal {E}_{N,1}^{int}(\cdot)$ is defined by
\begin{eqnarray*}
\mathcal {E}_{N,1}^{int}(\bar{V})=\sum_{|\alpha|\leq N-1}\langle
\partial^{\alpha}\bar{v},\nabla\partial^{\alpha}\bar{\sigma}\rangle.
\end{eqnarray*}
In fact, the first two equations of $ \eqref{sta.equ}$ can be
rewritten as
\begin{eqnarray}\label{3.7}
  &&\partial_t \bar{\sigma}+\nabla \cdot \bar{v}=f_{1},
\end{eqnarray}
\begin{eqnarray}\label{3.9}
  \partial_t \bar{v}+\nabla
  \bar{\sigma}+\frac{1}{\sqrt{\gamma}}\bar{E}=f_{2}-\frac{1}{\sqrt{\gamma}}\bar{v},
\end{eqnarray}
where
\begin{eqnarray}\label{f1f2}
 && \left\{
 \begin{aligned}
 & f_{1}:=-\bar{v}\cdot
  \nabla\bar{\sigma}-\frac{\gamma-1}{2}\bar{\sigma}\nabla \cdot \bar{v}
  -\bar{v}\cdot
 \nabla \sigma_{st}-\frac{\gamma-1}{2}\sigma_{st}\nabla \cdot \bar{v},\\
 & f_{2}:=-\bar{v}\cdot \nabla
  \bar{v}-\frac{\gamma-1}{2}\bar{\sigma}\nabla \bar{\sigma}-\bar{v}\times \bar{B}
  -\frac{\gamma-1}{2}\sigma_{st}\nabla \bar{\sigma}-\frac{\gamma-1}{2}\bar{\sigma}\nabla\sigma_{st}.
\end{aligned}\right.
\end{eqnarray}
Let $|\alpha|\leq N-1$. Applying $\partial^{\alpha}$ to (\ref{3.9}),
multiplying it by $\partial^{\alpha}\nabla \bar{\sigma}$, taking
integrations in $x$ and then using integration by parts and also the
final equation of (\ref{sta.equ}), replacing
$\partial_{t}\bar{\sigma}$ from (\ref{3.7}) give
\begin{eqnarray*}
\arraycolsep=1.5pt
\begin{array}{rl}
 &\displaystyle
\frac{d}{dt}\langle
\partial^{\alpha}\bar{v},\nabla
\partial^{\alpha}\bar{\sigma}\rangle
+\|\nabla\partial^{\alpha}\bar{\sigma}\|^{2}+\frac{1}{\gamma}\|
\partial^{\alpha}\bar{\sigma}\|^2\\[3mm]
=&\displaystyle -\frac{1}{\gamma}\langle
\partial^{\alpha}\left(\Phi(\bar{\sigma}+\sigma_{st})-\Phi(\sigma_{st})\right),
\partial^{\alpha}\bar{\sigma}\rangle+\langle\partial^{\alpha}f_{2},\nabla\partial^{\alpha}\bar{\sigma}\rangle\\[3mm]
&-\displaystyle\frac{1}{\sqrt{\gamma}}\langle\partial^{\alpha}\bar{v},\nabla\partial^{\alpha}\bar{\sigma}\rangle
+\|\nabla \cdot
\partial^{\alpha}\bar{v}\|^{2}-\langle\partial^{\alpha}f_{1},\nabla \cdot
\partial^{\alpha}\bar{v}\rangle.
\end{array}
\end{eqnarray*}
Then, it follows from Cauchy-Schwarz inequality that
\begin{eqnarray}\label{3.11}
\arraycolsep=1.5pt
\begin{array}{rl}
 &\displaystyle
\frac{d}{dt}\langle
\partial^{\alpha}\bar{v},\nabla
\partial^{\alpha}\bar{\sigma}\rangle
+\lambda(\|\nabla\partial^{\alpha}\bar{\sigma}\|^{2}+\|
\partial^{\alpha}\bar{\sigma}\|^2)\\[3mm]
\leq & C \displaystyle \|\nabla \cdot
\partial^{\alpha}\bar{v}\|^{2}+C(\|\partial^{\alpha}
\left(\Phi(\bar{\sigma}+\sigma_{st})-\Phi(\sigma_{st})\right) \|^{2}
+\|\partial^{\alpha}f_{1}\|^{2}+\|\partial^{\alpha}f_{2}\|^{2}).
\end{array}
\end{eqnarray}
Noticing that $\Phi(\sigma)$ is smooth in $\sigma$ with
$\Phi'(0)=0$, one has from (\ref{f1f2}) that
\begin{eqnarray*}
\arraycolsep=1.5pt
\begin{array}{rl}
&\|\partial^{\alpha}
\left(\Phi(\bar{\sigma}+\sigma_{st})-\Phi(\sigma_{st})\right) \|^{2}
+\|\partial^{\alpha}f_{1}\|^{2}+\|\partial^{\alpha}f_{2}\|^{2}\\[3mm]
\leq &C(\|[\bar{\sigma},\bar{v},\bar{B}]\|^{2}_{N}+\delta)\|
\nabla[\bar{\sigma},\bar{v}]\|_{N-1}^{2}.
\end{array}
\end{eqnarray*}
Here, if there is no derivative on $\bar{\sigma}$ or $\bar{v}$, then
use the $L^{\infty}$ of $\bar{\sigma}$ or $\bar{v}$. Plugging this
into (\ref{3.11}) taking summation over $|\alpha|\leq N-1$ yield
(\ref{step2}).

\medskip

\textbf{Step 3.} It holds that
\begin{equation}\label{step3}
  \begin{aligned}
  \dfrac{d}{dt}\mathcal {E}_{N,2}^{int}(\bar{V})+\lambda\|\bar{E}\|^{2}_{N-1} \leq
  C&\|[\bar{\sigma},\bar{v}]\|_{N}^{2}+C\|\bar{v}\|_{N}\|\nabla
  \bar{B}\|_{N-2}\\
  &+C(\|[\bar{\sigma},\bar{v},\bar{B}]\|_{N}^{2}+\delta)
  \|\nabla[\bar{\sigma},\bar{v}]\|_{N-1}^{2},
  \end{aligned}
\end{equation}
where $\mathcal {E}_{N,2}^{int}(\cdot)$ is defined by
\begin{eqnarray*}
\mathcal {E}_{N,2}^{int}(\bar{V})=\sum_{|\alpha|\leq N-1}\langle
\partial^{\alpha}\bar{v},\partial^{\alpha}\bar{E}\rangle.
\end{eqnarray*}
Applying $\partial^{\alpha}$ to (\ref{3.9}), multiplying it by
$\partial^{\alpha}\bar{E}$, taking integrations in $x$ and using
integration by parts and replacing $ \partial_{t}\bar{E}$  from the
third equation of (\ref{sta.equ}) give
\begin{equation*}
\arraycolsep=1.5pt
\begin{array}{rl}
 &\dfrac{d}{dt}\langle
\partial^{\alpha}\bar{v},\partial^{\alpha}\bar{E}\rangle+
\dfrac{1}{\sqrt{\gamma}}\|\partial^{\alpha}\bar{E}\|^{2}\\[3mm]
=& \dfrac{1}{\sqrt{\gamma}}\|\partial^{\alpha}\bar{v}\|^{2}+
\dfrac{1}{\sqrt{\gamma}}\langle
\partial^{\alpha}\bar{v},\nabla \times \partial^{\alpha}\bar{B}\rangle+\dfrac{1}{\sqrt{\gamma}}\langle
\partial^{\alpha}\bar{v},
\partial^{\alpha}[\Phi(\bar{\sigma}+\sigma_{st})\bar{v}+(\bar{\sigma}+\sigma_{st})\bar{v}]\rangle\\[3mm]
&-\langle\partial^{\alpha}\nabla
\bar{\sigma}+\dfrac{1}{\sqrt{\gamma}}\partial^{\alpha}\bar{v},\partial^{\alpha}\bar{E}\rangle
+\langle\partial^{\alpha}f_{2},
\partial^{\alpha}\bar{E}\rangle,
\end{array}
\end{equation*}
which from the Cauchy-Schwarz inequality further implies
\begin{equation*}
\arraycolsep=1.5pt
\begin{array}{rl}
 &\dfrac{d}{dt}\langle
\partial^{\alpha}\bar{v},\partial^{\alpha}\bar{E}\rangle+
\lambda\|\partial^{\alpha}\bar{E}\|^{2}\\[3mm]
 \leq &
  C\|[\bar{\sigma},\bar{v}]\|_{N}^{2}+C\|\bar{v}\|_{N}\|\nabla
  \bar{B}\|_{N-2}+C(\|[\bar{\sigma},\bar{v},\bar{B}]\|_{N}^{2}+\delta)
  \|\nabla[\bar{\sigma},\bar{v}]\|_{N-1}^{2}.
\end{array}
\end{equation*}
Thus $\eqref{step3}$ follows from taking summation of the above
estimate over $|\alpha|\leq N-1$.

\medskip

\textbf{Step 4.} It holds that
\begin{equation}\label{step4}
\begin{aligned}
   \frac{d}{dt}\mathcal {E}_{N,3}^{int}(\bar{V})+\lambda\|\nabla\bar{B}\|^{2}_{N-2}
   \leq & C\|[\bar{v},\bar{E}]\|_{N-1}^{2}\\
   &+C(\|\bar{\sigma}\|_{N}^{2}+\delta)\|\nabla \bar{v}\|_{N-1}^{2},
\end{aligned}
\end{equation}
where $\mathcal {E}_{N,3}^{int}(\cdot)$ is defined by
\begin{eqnarray*}
\mathcal {E}_{N,3}^{int}(\bar{V})=-\sum_{|\alpha|\leq N-2}\langle
\nabla \times
\partial^{\alpha}\bar{E},\partial^{\alpha}\bar{B}\rangle.
\end{eqnarray*}
In fact, for $|\alpha|\leq N-2$, applying $\partial^{\alpha}$ to the
third equation of  $\eqref{sta.equ}$, multiplying it by
$-\partial^{\alpha}\nabla \times \bar{B}$, taking integrations in
$x$ and using integration by parts and replacing $
\partial_{t}\bar{B}$  from the fourth  equation of $\eqref{sta.equ}$
implie
\begin{equation*}
\arraycolsep=1.5pt
\begin{array}{rl}
& -\dfrac{d}{dt}\langle
\partial^{\alpha}\bar{E},\nabla \times \partial^{\alpha}\bar{B}\rangle+
\dfrac{1}{\sqrt{\gamma}}\|\nabla\times \partial^{\alpha}\bar{B}\|^{2} \\[3mm]
=&\dfrac{1}{\sqrt{\gamma}}\|\nabla\times\partial^{\alpha}\bar{E}\|^{2}-\dfrac{1}{\sqrt{\gamma}}
\langle \partial^{\alpha}\bar{v},\nabla \times
\partial^{\alpha}\bar{B}\rangle
-\dfrac{1}{\sqrt{\gamma}}\langle
\partial^{\alpha}[\Phi(\bar{\sigma}+\sigma_{st})\bar{v}+(\bar{\sigma}+\sigma_{st})\bar{v}],\nabla \times
\partial^{\alpha}\bar{B}\rangle,
\end{array}
\end{equation*}
which gives $\eqref{step4}$ by further using Cauchy-Schwarz
inequality and taking summation over $|\alpha|\leq N-2$, where we
also used
\begin{eqnarray*}
\|\partial^{\alpha}\partial_{i}\bar{B}\|=\|\partial_{i}\Delta^{-1}\nabla
\times(\nabla\times\partial^{\alpha}\bar{B}) \|\leq\|\nabla\times
\partial^{\alpha}\bar{B}\|
\end{eqnarray*}
for each $1\leq i\leq 3$, due to the fact
$\partial_{i}\Delta^{-1}\nabla$ is bounded from $L^{p}$ to itself
for $1<p<\infty$, cf. \cite{Stein}.

\medskip

 \textbf{Step 5.} Now, following the four steps above, we are ready to prove
 $\eqref{3.2}$. Let us define
\begin{eqnarray*}
   \mathcal {E}_{N}(\bar{V}(t))=\sum_{|\alpha|\leq
N}\int_{\mathbb{R}^3}(1+\sigma_{st}+\Phi(\sigma_{st}))
(|\partial^{\alpha}\bar{\sigma}|^2+|\partial^{\alpha}\bar{v}|^2)dx+\|[\bar{E},\bar{B}]\|_{N}^{2}
   +\sum_{i=1}^{3}\kappa_{i}\mathcal {E}^{int}_{N,i}(\bar{V}(t)),
\end{eqnarray*}
that is,
\begin{equation}\label{3.12}
\arraycolsep=1.5pt
\begin{array}{rl}
\mathcal{E}_{N}(\bar{V}(t))=&\displaystyle\sum_{|\alpha|\leq
N}\int_{\mathbb{R}^3}(1+\sigma_{st}+\Phi(\sigma_{st}))
(|\partial^{\alpha}\bar{\sigma}|^2+|\partial^{\alpha}\bar{v}|^2)dx+\|[\bar{E},\bar{B}]\|_{N}^{2}\\[3mm]
&\displaystyle+\kappa_{1}\sum_{|\alpha|\leq N-1} \langle
\partial^{\alpha}\bar{v},\nabla\partial^{\alpha}\bar{\sigma}\rangle+\kappa_{2}\sum_{|\alpha|\leq N-1}\langle
\partial^{\alpha}\bar{v},\partial^{\alpha}\bar{E}\rangle\\[3mm]
&\displaystyle-\kappa_{3}\sum_{|\alpha|\leq N-2}\langle \nabla
\times
\partial^{\alpha}\bar{E},\partial^{\alpha}\bar{B}\rangle
\end{array}
\end{equation}
for constants $0<\kappa_{3}\ll\kappa_{2}\ll\kappa_{1}\ll 1$ to be
determined. Notice that as long as $ 0<\kappa_{i}\ll 1$ is small
enough for $i=1,2,3$, and $\sigma_{st}+\Phi(\sigma_{st}) $ depending
only on $ x$ is sufficiently small compared with $1$, then
$\mathcal{E}_{N}(\bar{V}(t))\sim \|\bar{V}(t)\|^{2}_{N}$ holds true.
Moreover, letting $0<\kappa_{3}\ll\kappa_{2}\ll\kappa_{1}\ll 1$ with
$\kappa_{2}^{3/2}\ll\kappa_{3}$, the sum of $\eqref{3.3}\times
\kappa_{1}$, $\eqref{step2}\times \kappa_{2}$, $\eqref{step4}\times
\kappa_{3}$ implies that there are $\lambda>0$, $C>0$ such that
$\eqref{3.2}$ holds true with $\mathcal {D}_{N}(\cdot)$ defined in
$\eqref{de.D}$. Here, we have used the following Cauchy-Schwarz
inequality:
\begin{eqnarray*}
   2 \kappa_{2} \|\bar{v}\|_{N}\|\nabla \bar {B}\|_{N-2}\leq
   \kappa_{2}^{1/2}\|\bar{v}\|_{N}^{2}+\kappa_{2}^{3/2}\|\nabla\bar{B}\|^{2}_{N-2}.
\end{eqnarray*}
Due to $\kappa_{2}^{3/2}\ll \kappa_{3}$, both terms on the r.h.s. of
the above inequality were absorbed. This completes the proof of
Theorem $\ref{estimate}$.
\end{proof}

Since $\eqref{sta.equ}$ is a quasi-linear symmetric hyperbolic
system, the short-time existence can be proved in much more general
case as in \cite{Kato}; see also (Theorem 1.2, Proposition 1.3, and
Proposition 1.4 in Chapter 16 of \cite{Taylor}). From Theorem
\ref{estimate} and the continuity argument, it is easy to see that $
\mathcal {E}_{N}(\bar{V}(t)) $ is bounded uniformly in time under
the assumptions that $\mathcal {E}_{N}(\bar{V}_{0})>0$ and
$\|n_{b}-1\|_{W_{0}^{N+1,2}}$ are small enough. Therefore, the
global existence of solutions satisfying \eqref{V.satisfy} and
\eqref{pro.2.1j} follows in the standard way; see also \cite{Duan}.
This completes the proof of Proposition \ref{pro.2.1}.\qed

\vspace{5mm}

\section{Decay in time for the non-linear system}\label{sec4}
In this section, we are devoted to the rate of the convergence of
solution to the equilibrium $[n_{st},0,E_{st},0]$ for the system
\eqref{1.1} over $\mathbb{R}^3$. In fact by setting
\begin{eqnarray*}
\bar{\rho}=n-n_{st},\ \ \bar{u}=u,\ \ E_{1}=E-E_{st},\ \ B_{1}=B,
\end{eqnarray*}
and
\begin{eqnarray*}
\rho_{st}=n_{st}-1,
\end{eqnarray*}
then $\bar{U}:=[\bar{\rho},\bar{u},E_{1},B_{1}]$ satisfies
\begin{equation}\label{rhost}
\left\{
  \begin{aligned}
  &\partial_t \bar{\rho}+\nabla\cdot \bar{u}=g_{1} ,\\
  &\partial_t \bar{u}+\bar{u} + E_{1} +\gamma \nabla\bar{\rho}=g_{2},
  \\
  &\partial_t E_{1}-\nabla\times B_{1}-\bar{u}=g_{3},\\
  &\partial_t B_{1}+\nabla \times E_{1}=0,\\
  &\nabla \cdot E_{1}=-\bar{\rho}, \ \  \nabla \cdot B_{1}=0, \ \ \ t>0,\ x\in\mathbb{R}^{3},\\
\end{aligned}\right.
\end{equation}
with initial data
\begin{eqnarray}\label{rhosti}
\begin{aligned}
\bar{U}|_{t=0}=\bar{U}_{0}:=&[\bar{\rho}_{0},\bar{u}_{0},E_{1,0},B_{1,0}]\\
 =&[n_0-n_{st},u_0,E_{0}-E_{st},B_0], \ \
\ x\in\mathbb{R}^{3},
\end{aligned}
\end{eqnarray}
satisfying the compatible conditions
\begin{eqnarray}\label{rhostC}
\nabla \cdot E_{1,0}=-\bar{\rho}_{0}, \ \  \nabla \cdot B_{1,0}=0.
\end{eqnarray}
Here the nonlinear source terms take the form of
\begin{equation}\label{sec5.ggg}
\arraycolsep=1.5pt \left\{
 \begin{aligned}
 & g_{1}=-\nabla\cdot[(\bar{\rho}+\rho_{st}) \bar{u}],\\
 &\begin{array}[b]{rcl}
  g_{2}&=&-\bar{u} \cdot \nabla \bar{u}-\bar{u}\times B_{1}
  -\gamma [(\bar{\rho}+1+\rho_{st})^{\gamma-2}-1]\nabla\bar{\rho}\\
  &&-\gamma
  [(1+\bar{\rho}+\rho_{st})^{\gamma-2}-(1+\rho_{st})^{\gamma-2}]\nabla\rho_{st},
  \end{array}\\
 & g_{3}=(\bar{\rho}+\rho_{st}) \bar{u}.
\end{aligned}\right.
\end{equation}

In what follows, we will denote $[\rho,u,E,B]$ as the solution to
the  the following linearized equation of \eqref{rhost}:
\begin{equation}\label{DJ}
\left\{
  \begin{aligned}
  &\partial_t \rho+\nabla\cdot u=0,\\
  &\partial_t u+u+ E +\gamma \nabla\rho=0,\\
  &\partial_t E-\nabla\times B-u=0,\\
  &\partial_t B+\nabla \times E=0,\\
  &\nabla \cdot E=-\rho, \ \  \nabla \cdot B=0, \ \ \ t>0, \ \ x\in\mathbb{R}^{3},\\
\end{aligned}\right.
\end{equation}
with given initial data
\begin{eqnarray}\label{2.61}
U|_{t=0}=\bar{U}_{0}:=[\bar{\rho}_{0},\bar{u}_{0},E_{1,0},B_{1,0}],
\ \ \ x\in\mathbb{R}^{3},
\end{eqnarray}
satisfying the compatible conditions \eqref{rhostC}.

For the above linearized equations, the $L^{p}$-$L^{q}$ time-decay
property was proved by Duan in \cite{Duan}. We  list only some
special $L^{p}$-$L^{q}$ time decay properties  in the following
proposition.
\begin{proposition}\label{thm.decay}
 Suppose $U(t)=e^{tL}\bar{U}_{0}$ is the solution
to the Cauchy problem \eqref{DJ}-\eqref{2.61} with the initial data
$\bar{U}_{0}=[\bar{\rho}_{0},\bar{u}_{0},E_{1,0},B_{1,0}] $
satisfying \eqref{rhostC}. Then, $U=[\rho,u,E,B]$ satisfies the
following time-decay property:

\begin{eqnarray}\label{col.decay1}
 && \left\{
 \begin{aligned}
 & \|\rho(t)\|\leq C e^{-\frac{t}{2}}\|[\bar{\rho}_{0},\bar{u}_{0}]\|,\\
 & \|u(t)\| \leq C e^{-\frac{t}{2}}\|\bar{\rho}_{0}\|+C(1+t)^{-\frac{5}{4}}
 \|[\bar{u}_{0}, E_{1,0},B_{1,0}]\|_{L^1\cap \dot{H}^{2}},\\
 &\|E(t)\|\leq C (1+t)^{-\frac{5}{4}}
 \|[\bar{u}_{0}, E_{1,0},B_{1,0}]\|_{L^1\cap \dot{H}^{3}},\\
 &\|B(t)\|\leq C (1+t)^{-\frac{3}{4}}
 \|[\bar{u}_{0}, E_{1,0},B_{1,0}]\|_{L^1\cap \dot{H}^{2}},
\end{aligned}\right.
\end{eqnarray}
and
\begin{eqnarray}\label{col.decayinfty1}
 && \left\{
 \begin{aligned}
 & \|\rho(t)\|_{\infty}\leq C e^{-\frac{t}{2}}\|[\bar{\rho}_{0},\bar{u}_{0}]\|_{L^{2}\cap\dot{H}^{2}},\\
 & \|u(t)\|_{\infty} \leq C e^{-\frac{t}{2}}\|\bar{\rho}_{0}\|_{L^{2}\cap\dot{H}^{2}}+C(1+t)^{-2}
 \|[\bar{u}_{0}, E_{1,0},B_{1,0}]\|_{L^1\cap \dot{H}^{5}},\\
 &\|E(t)\|_{\infty}\leq C (1+t)^{-2}
 \|[\bar{u}_{0}, E_{1,0},B_{1,0}]\|_{L^1\cap \dot{H}^{6}},\\
 &\|B(t)\|_{\infty}\leq C (1+t)^{-\frac{3}{2}}
 \|[\bar{u}_{0}, E_{1,0},B_{1,0}]\|_{L^1\cap \dot{H}^{5}},
\end{aligned}\right.
\end{eqnarray}
and, moreover,
\begin{eqnarray}\label{col.EB}
 && \left\{
 \begin{aligned}
 &\|\nabla B(t)\|\leq C (1+t)^{-\frac{5}{4}}
 \|[\bar{u}_{0}, E_{1,0},B_{1,0}]\|_{ L^1 \cap \dot{H}^{4}},\\
 & \|\nabla^{N}[E(t),B(t)]\|\leq C(1+t)^{-\frac{5}{4}}
 \|[\bar{u}_{0}, E_{1,0},B_{1,0}]\|_{ L^1 \cap \dot{H}^{N+3}}.
\end{aligned}\right.
\end{eqnarray}

\end{proposition}

In what follows, since we shall apply the linear $L^{p}$-$L^{q}$
time-decay property of the homogeneous system \eqref{DJ}, we need
the mild form of the non-linear Cauchy problem
\eqref{rhost}-\eqref{rhosti}. From now on, we always denote
$\bar{U}=[\bar{\rho},\bar{u},E_{1},B_{1}]$ to the non-linear Cauchy
problem $\eqref{rhost}$-$\eqref{rhosti}$. Then, by Duhamel's
principle, the solution $\bar{U}$ can be formally written as
\begin{eqnarray}\label{sec5.U}
\bar{U}(t)=e^{tL}\bar{U}_{0}+\int_{0}^{t}e^{(t-s)L}[g_{1}(s),g_{2}(s),g_{3}(s),0]d
s,
\end{eqnarray}
where $e^{tL}\bar{U}_{0}$ denotes the solution to the Cauchy problem
$\eqref{DJ}$-$\eqref{2.61}$ without nonlinear sources.

The following  two lemmas give the full and high-order energy
estimates.
\begin{lemma}\label{lem.V}
Let $\bar{V}=[\bar{\sigma},\bar{v},\bar{E},\bar{B}]$ be the solution
to the Cauchy problem $\eqref{sta.equ}$--$ \eqref{sta.equi}$ with
initial data
$\bar{V}_{0}=[\bar{\sigma}_{0},\bar{v}_{0},\bar{E}_{0},\bar{B}_{0}]$
satisfying $\eqref{sta.equC}$. Then, if $\mathcal
{E}_{N}(\bar{V}_{0})$ and $\|n_{b}-1\|_{W_{0}^{N+1,2}}$ are
sufficiently small,
\begin{eqnarray}\label{sec5.ENV0}
\dfrac{d}{dt}\mathcal {E}_{N}(\bar{V}(t))+\lambda \mathcal {D}_{N}(\bar{V}(t))
\leq 0
\end{eqnarray}
holds for any $t>0$, where $\mathcal {E}_{N}(\bar{V}(t))$, $\mathcal
{D}_{N}(\bar{V}(t))$ are defined in the form of $\eqref{de.E}$ and
$\eqref{de.D}$, respectively.
\end{lemma}
\begin{proof}
It can be seen directly from the proof of Theorem \ref{estimate}.
\end{proof}

\begin{lemma}\label{estimate2}
Let $\bar{V}=[\bar{\sigma},\bar{v},\bar{E},\bar{B}]$ be the solution
to the Cauchy problem $\eqref{sta.equ}$-$\eqref{sta.equi}$ with
initial data
$\bar{V}_{0}=[\bar{\sigma}_{0},\bar{v}_{0},\bar{E}_{0},\bar{B}_{0}]$
satisfying $\eqref{sta.equC}$ in the sense of Proposition
$\ref{pro.2.1}$. Then if $ \mathcal {E}_{N}(\bar{V}_{0})$ and
$\|n_{b}-1\|_{W_{0}^{N+1,2}}$ are sufficiently small, there are the
high-order instant energy functional $\mathcal {E}_{N}^{h}(\cdot)$
and the corresponding dissipation rate $\mathcal {D}_{N}^{h}(\cdot)$
such that
\begin{eqnarray}\label{sec5.high}
&& \frac{d}{dt}\mathcal {E}_{N}^{h}(\bar{V}(t))+\lambda\mathcal
{D}^{h}_{N}(\bar{V}(t))\leq 0,
\end{eqnarray}
holds for any $ t \geq 0$.
\end{lemma}

\begin{proof}
The proof can be done by modifying the proof of Theorem
$\ref{estimate}$ a little. In fact, by letting the energy estimates
made only on the high-order derivatives, then corresponding to
$\eqref{3.3}$, $\eqref{step2}$, $\eqref{step3}$ and $\eqref{step4}$,
it can be re-verified that

\begin{equation*}
\arraycolsep=1.5pt
\begin{array}{rl}
&\displaystyle \frac{1}{2}\frac{d}{dt}\left(\sum_{1\leq|\alpha|\leq
N}\int_{\mathbb{R}^3}(1+\sigma_{st}+\Phi(\sigma_{st}))
(|\partial^{\alpha}\bar{\sigma}|^2+|\partial^{\alpha}\bar{v}|^2)dx+\|\nabla[\bar{E},\bar{B}]\|_{N-1}^{2}\right)\\[5mm]
&\displaystyle+\frac{1}{\sqrt{\gamma}}\sum_{1\leq|\alpha|\leq
N}\int_{\mathbb{R}^3}(1+\sigma_{st}+\Phi(\sigma_{st}))|\partial^{\alpha}\bar{v}|^{2}dx\\[5mm]
\leq & \displaystyle
C(\|\bar{V}\|_{N}+\delta)(\|\nabla[\bar{\sigma},\bar{v}]\|_{N-1}^{2}
 +\|\nabla \bar{E}\|_{N-2}^2),
 \end{array}
\end{equation*}
\begin{eqnarray*}
 \frac{d}{dt}\sum_{1\leq|\alpha|\leq N-1}\langle
\partial^{\alpha}\bar{v},\nabla\partial^{\alpha}\bar{\sigma}\rangle+\lambda\|\nabla\bar{\sigma}\|^{2}_{N-1}
  \leq  C\|\nabla^2\bar{v}\|_{N-2}^{2}+C(\|[\bar{\sigma}, \bar{v},\bar{B}]\|_{N}^{2}+\delta)
  \|\nabla[\bar{\sigma},\bar{v}]\|_{N-1}^{2},
\end{eqnarray*}
\begin{eqnarray*}
\begin{aligned}
  \dfrac{d}{dt}\sum_{1\leq|\alpha|\leq N-1}\langle
\partial^{\alpha}\bar{v},\partial^{\alpha}\bar{E}\rangle+\lambda\|\nabla\bar{E}\|^{2}_{N-2} \leq
  C&\|\nabla[\bar{\sigma},\bar{v}]\|_{N-1}^{2}+C\|\nabla\bar{v}\|_{N-1}\|\nabla^2
  \bar{B}\|_{N-3}\\
  &+C(\|[\bar{\sigma},\bar{v},\bar{B}]\|_{N}^{2}+\delta)
  \|\nabla[\bar{\sigma},\bar{v}]\|_{N-1}^{2},
  \end{aligned}
\end{eqnarray*}
and
\begin{eqnarray*}
&&\begin{aligned}
  & -\frac{d}{dt}\sum_{1\leq |\alpha|\leq N-2}\langle \nabla
   \times\partial^{\alpha}\bar{E},\partial^{\alpha}\bar{B}\rangle+\lambda\|\nabla^{2}\bar{B}\|^{2}_{N-3}\\
    \leq & C\|\nabla^2\bar{E}\|_{N-3}^{2}
   +C\|\nabla\bar{v}\|_{N-3}^2+C(\|\bar{\sigma}\|_{N}^{2}+\delta)\|\nabla \bar{v}\|_{N-1}^{2}.
\end{aligned}
\end{eqnarray*}
Here, the details of proof are omitted for simplicity. Now, similar
to $\eqref{3.12}$, let us define
\begin{equation}\label{def.high}
\begin{aligned}
\mathcal{E}_{N}^{h}(\bar{V}(t))&=\sum_{1\leq|\alpha|\leq
N}\int_{\mathbb{R}^3}(1+\sigma_{st}+\Phi(\sigma_{st}))
(|\partial^{\alpha}\bar{\sigma}|^2+|\partial^{\alpha}\bar{v}|^2)dx+\|\nabla[\bar{E},\bar{B}]\|_{N-1}^{2}\\
&+\kappa_{1}\sum_{1\leq|\alpha|\leq N-1}\langle
\partial^{\alpha}\bar{v},\nabla\partial^{\alpha}\bar{\sigma}\rangle+\kappa_{2}\sum_{1\leq|\alpha|\leq N-1}\langle
\partial^{\alpha}\bar{v},\partial^{\alpha}\bar{E}\rangle\\[3mm]
&-\kappa_{3}\sum_{1\leq |\alpha|\leq N-2}\langle \nabla
   \times\partial^{\alpha}\bar{E},\partial^{\alpha}\bar{B}\rangle.
\end{aligned}
\end{equation}

Similarly, one can choose $0<\kappa_{3}\ll\kappa_{2}\ll\kappa_{1}\ll
1$ with $\kappa_{2}^{3/2}\ll\kappa_{3}$ such that $\mathcal
{E}_{N}^{h}(\bar{V}(t))\sim \|\nabla \bar{V}(t)\|_{N-1}^{2}$ because
$\sigma_{st}+\Phi(\sigma_{st}) $ depends only on $ x$ sufficiently
small compared with $1$. Furthermore, the linear combination of
previously obtained four estimates with coefficients corresponding
to $\eqref{def.high}$ yields $\eqref{sec5.high}$ with $\mathcal
{D}_{N}^{h}(\cdot)$ defined in $\eqref{de.Dh}$. This completes the
proof of Lemma \ref{estimate2}.
\end{proof}

Now, we begin with the time-weighted estimate and iteration for the
Lyapunov inequality $\eqref{sec5.ENV0}$. Let $\ell \geq 0$.
Multiplying $\eqref{sec5.ENV0}$ by $(1+t)^{\ell}$ and taking
integration over $[0,t]$ give

\begin{eqnarray*}
\begin{aligned}
 & (1+t)^{\ell}\mathcal {E}_{N}(\bar{V}(t))+\lambda
 \int_{0}^{t}(1+s)^{\ell}\mathcal {D}_{N}(\bar{V}(s))d s \\
 \leq & \mathcal {E}_{N}(\bar{V}_{0})+ \ell
 \int_{0}^{t}(1+s)^{\ell-1}\mathcal {E}_{N}(\bar{V}(s))d s.
\end{aligned}
\end{eqnarray*}
Noticing
\begin{eqnarray*}
  \mathcal {E}_{N}(\bar{V}(t))
  \leq C (D_{N+1}(\bar{V}(t))+\|
  \bar{B}\|^{2}),
\end{eqnarray*}
it follows that
\begin{eqnarray*}
\begin{aligned}
 & (1+t)^{\ell}\mathcal {E}_{N}(\bar{V}(t))+\lambda
 \int_{0}^{t}(1+s)^{\ell}\mathcal {D}_{N}(\bar{V}(s))d s \\
 \leq & \mathcal {E}_{N}(\bar{V}_{0})+ C \ell
 \int_{0}^{t}(1+s)^{\ell-1}\|
  \bar{B}(s)\|^{2}d s+ C\ell\int_{0}^{t}(1+s)^{\ell-1}\mathcal {D}_{N+1}(\bar{V}(s))d
  s.
\end{aligned}
\end{eqnarray*}
Similarly, it holds that
\begin{eqnarray*}
\begin{aligned}
 & (1+t)^{\ell-1}\mathcal {E}_{N+1}(\bar{V}(t))+\lambda
 \int_{0}^{t}(1+s)^{\ell-1}\mathcal {D}_{N+1}(\bar{V}(s))d s \\
 \leq & \mathcal {E}_{N+1}(\bar{V}_{0})+ C (\ell-1)
 \int_{0}^{t}(1+s)^{\ell-2}\|
  \bar{B}(s)\|^{2}ds
+ C(\ell-1)\int_{0}^{t}(1+s)^{\ell-2}\mathcal {D}_{N+2}(\bar{V}(s))d s,
\end{aligned}
\end{eqnarray*}
and
\begin{eqnarray*}
\mathcal {E}_{N+2}(\bar{V}(t))+\lambda \int_{0}^{t}\mathcal
{D}_{N+2}(\bar{V}(s))d s \leq \mathcal {E}_{N+2}(\bar{V}_{0}).
\end{eqnarray*}
Then, for $1<\ell<2$, it follows by iterating the above estimates
that
\begin{eqnarray}\label{sec5.ED}
\begin{aligned}
 & (1+t)^{\ell}\mathcal {E}_{N}(\bar{V}(t))+\lambda
 \int_{0}^{t}(1+s)^{\ell}\mathcal {D}_{N}(\bar{V}(s))d s \\
 \leq & C \mathcal {E}_{N+2}(\bar{V}_{0})+ C
 \int_{0}^{t}(1+s)^{\ell-1}\|
  \bar{B}(s)\|^{2}d s.
\end{aligned}
\end{eqnarray}

Similarly, for $2<\kappa<3$, the time-weighted estimate and
iteration for the Lyapunov inequality $\eqref{sec5.high}$ give
\begin{eqnarray*}
\begin{aligned}
 & (1+t)^{\kappa}\mathcal {E}_{N}^h(\bar{V}(t))+\lambda
 \int_{0}^{t}(1+s)^{\kappa}\mathcal {D}_{N}^h(\bar{V}(s))d s \\
 \leq & C \mathcal {E}_{N+3}^h(\bar{V}_{0})+ C
 \int_{0}^{t}(1+s)^{\kappa-1}\|
  \nabla\bar{B}(s)\|^{2}d s.
\end{aligned}
\end{eqnarray*}
Here the smallness of $\|n_{b}-1\|_{W_{0}^{N+4,2}}$ has been used in
the process of iteration for the Lyapunov inequalities
$\eqref{sec5.ENV0}$ and $\eqref{sec5.high}$. Taking $\kappa=l+1$, it
holds that
\begin{multline}\label{sec5.EhD}
(1+t)^{l+1}\mathcal {E}_{N}^h(\bar{V}(t))+\lambda
 \int_{0}^{t}(1+s)^{l+1}\mathcal {D}_{N}^h(\bar{V}(s))d s \\
 \leq  C \mathcal {E}_{N+3}^h(\bar{V}_{0})+ C
 \int_{0}^{t}(1+s)^{l}\|
  \nabla\bar{B}(s)\|^{2}d s\\
   \leq C \mathcal {E}_{N+3}^h(\bar{V}_{0})+ C \int_{0}^{t}(1+s)^{\ell}\mathcal {D}_{N}(\bar{V}(s))d s.
\end{multline}
Combining $\eqref{sec5.ED}$ with $\eqref{sec5.EhD}$, we have
\begin{multline}\label{sec5.EDEhD}
(1+t)^{\ell}\mathcal {E}_{N}(\bar{V}(t))+
 \int_{0}^{t}(1+s)^{\ell}\mathcal {D}_{N}(\bar{V}(s))d s\\
 +(1+t)^{l+1}\mathcal
{E}_{N}^h(\bar{V}(t))+
 \int_{0}^{t}(1+s)^{l+1}\mathcal {D}_{N}^h(\bar{V}(s))d s \\
 \leq  C \mathcal {E}_{N+3}(\bar{V}_{0})+ C \int_{0}^{t}(1+s)^{\ell-1}\|
  \bar{B}(s)\|^{2}d s.
\end{multline}

 For this time, to estimate the integral term on the r.h.s. of
$\eqref{sec5.EDEhD}$, let's define
\begin{eqnarray}\label{sec5.def}
\mathcal {E}_{N,\infty}(\bar{V}(t))=\sup\limits_{0\leq s \leq t} \
\left\{(1+s)^{\frac{3}{2}}\mathcal
{E}_{N}(\bar{V}(s))+(1+s)^{\frac{5}{2}}\mathcal
{E}_{N}^h(\bar{V}(s))\right\},
\end{eqnarray}
\begin{eqnarray}\label{sec5.defL}
L_{0}(t)=\sup\limits_{0\leq s \leq t}
(1+s)^{\frac{5}{2}}\|[\bar{\rho},\bar{u}]\|^{2}.
\end{eqnarray}
Then, we have the following
\begin{lemma}\label{lem.Bsigma}
For any $t\geq0$, it holds that:
\begin{eqnarray}\label{lem.tildeB}
&&\begin{aligned} \|\bar{B}(t)\|^2\leq C
(1+t)^{-\frac{3}{2}}\left(\|[\bar{\sigma}_{0},\bar{v}_{0}]\|^{2}+
 \|[\bar{v}_{0},\right.& \bar{E}_{0},\bar{B}_{0}]\|^2_{L^1\cap \dot{H}^{2}}\\
&\left.+[\mathcal {E}_{N,\infty}(\bar{V}(t))]^2+\delta^2 \mathcal
{E}_{N,\infty}(\bar{V}(t))\right).
\end{aligned}
\end{eqnarray}
\end{lemma}

\begin{proof}
Applying the  fourth linear estimate on $B$ in $\eqref{col.decay1}$
to the mild form \eqref{sec5.U} gives
\begin{eqnarray}\label{sec5.decayB}
&&\begin{aligned} \|B_{1}(t)\|\leq C (1+t)^{-\frac{3}{4}}
 \|[\bar{u}_{0},& E_{1,0},B_{1,0}]\|_{L^1\cap \dot{H}^{2}}\\
&+C
\int_{0}^{t}(1+t-s)^{-\frac{3}{4}}\|[g_{2}(s),g_{3}(s)]\|_{L^{1}\cap\dot{H}^{2}}ds.
\end{aligned}
\end{eqnarray}
Applying the $L^{2}$ linear estimate on $u$ in $\eqref{col.decay1}$
to the mild form $\eqref{sec5.U}$,
\begin{eqnarray}\label{baruL2}
 &&\begin{aligned}
\|\bar{u}(t)\| \leq C(1+t)^{-\frac{5}{4}}(
 &\|\bar{\rho}_{0}\|+\|[\bar{u}_{0}, E_{1,0},B_{1,0}]\|_{L^{1}\cap\dot{H}^{2}})\\
  &+C \int_{0}^{t}(1+t-s)^{-\frac{5}{4}}\left(\|g_{1}(s)\|+\|[g_{2}(s),g_{3}(s)]\|_{L^{1}\cap
  \dot{H}^{2}}\right)ds.
 \end{aligned}
\end{eqnarray}
Applying the $L^{2}$ linear estimate on $\rho$ in
$\eqref{col.decay1}$ to $\eqref{sec5.U}$, one has
\begin{eqnarray}\label{rhoL2}
\|\bar{\rho}(t)\|\leq C
e^{-\frac{t}{2}}\|[\bar{\rho}_{0},\bar{u}_{0}]\|+ C
 \int_{0}^{t}e^{-\frac{t-s}{2}}\|[g_{1}(s),g_{2}(s)]\|d s.
\end{eqnarray}

Recall the definition $\eqref{sec5.ggg}$ of $g_{1}$, $g_{2}$ and
$g_{3}$,
\begin{eqnarray*}
\begin{aligned}
& g_{1}(s)=-\rho_{st} \nabla \cdot \bar{u}-\bar{\rho} \nabla \cdot
\bar{u}-\bar{u} \cdot \nabla \rho_{st}-\bar{u} \cdot \nabla
\bar{\rho},\\
& g_{2}(s)\sim \bar{u} \cdot \nabla \bar{u} + \bar{u}\times B_{1} +\bar{\rho}\nabla \bar{\rho}+ \rho_{st}\nabla \bar{\rho} +\bar{\rho}\nabla \rho_{st},\\
& g_{3}(s)=\bar{\rho} \bar{u}+\rho_{st} \bar{u}.
\end{aligned}
\end{eqnarray*}
Firstly, we estimate those terms including $\rho_{st}$. It follows that
\begin{eqnarray*}
\begin{aligned}
&\|\rho_{st}\nabla \cdot \bar{u}\|\leq
\|\rho_{st}\|_{L^{\infty}}\|\nabla \bar{u}\|,\ \ \ \ \ \|\bar{u}
\cdot \nabla \rho_{st}\|\leq
\|\nabla\rho_{st}\|\|\bar{u}\|_{L^{\infty}}\leq
\|\nabla\rho_{st}\|\|\nabla\bar{u}\|_{H^{1}},\\
&\|\rho_{st}\nabla \bar{\rho}\|_{L^1}\leq \|\rho_{st}\|\|\nabla
\bar{\rho}\|,\ \ \ \ \ \|\rho_{st}\nabla \bar{\rho}\|\leq
\|\rho_{st}\|_{L^{\infty}}\|\nabla\bar{\rho}\|\leq
\|\nabla\rho_{st}\|_{H^{1}}\|\nabla\bar{\rho}\|,\\
 & \|\bar{\rho}
\nabla \rho_{st}\|_{L^1}\leq \|\nabla\rho_{st}\|\|\bar{\rho}\|,\ \ \
\ \|\nabla\rho_{st}\bar{\rho}\|\leq
\|\bar{\rho}\|_{L^{\infty}}\|\rho_{st}\|\leq
\|\rho_{st}\|\|\nabla\bar{\rho}\|_{H^{1}},\\
&\|\rho_{st} \bar{u}\|_{L^1} \leq \|\rho_{st}\|\|\bar{u}\|,\\
\end{aligned}
\end{eqnarray*}
and for $|\alpha|=2$, one has
\begin{eqnarray*}
&&\begin{aligned} \|\partial^{\alpha}(\rho_{st}\nabla
\bar{\rho})\|\leq  & \|\rho_{st}\partial^{\alpha}\nabla\bar{\rho}\|
+\|\partial^{\alpha}(\rho_{st}\nabla
\bar{\rho})-\rho_{st}\partial^{\alpha}\nabla\bar{\rho}\|\\[3mm]
\leq &
\|\rho_{st}\|_{L^{\infty}}\|\partial^{\alpha}\nabla\bar{\rho}\|
+C\|\nabla\rho_{st}\|_{H^{|\alpha|-1}}\|\nabla\bar{\rho}\|_{L^{\infty}}+C\|\nabla
\rho_{st}\|_{L^{\infty}}\|\nabla \bar{\rho}\|_{H^{|\alpha|-1}}\\[3mm]
\leq & C\delta\|\nabla \bar{\rho}\|_{H^{2}},
\end{aligned}
\end{eqnarray*}
where we have used the estimate $\|\partial^{\alpha}(f g)-f
\partial^{\alpha}g\|\leq C\|\nabla f\|_{H^{k-1}}\|g\|_{L^{\infty}}+C\|\nabla f\|_{L^{\infty}}\|g\|_{H^{k-1}}$, for any $|\alpha|=k$.
 Similarly, it holds that
\begin{eqnarray*}
\begin{aligned}
\|\partial^{\alpha}(\rho_{st} \bar{u})\|\leq & \|\bar{u}
\partial^{\alpha}\rho_{st}\| +\|\partial^{\alpha}(\rho_{st}\bar{u})-\bar{u} \partial^{\alpha}\rho_{st}\|
\leq  C \delta \|\nabla \bar{u}\|_{H^{2}},
\end{aligned}
\end{eqnarray*}
\begin{eqnarray*}
\begin{aligned}
& \|\partial^{\alpha}(\bar{\rho} \nabla \rho_{st})\|\leq C \delta
\|\nabla \bar{\rho}\|_{H^{2}}.
\end{aligned}
\end{eqnarray*}
It is straightforward to verify that for  any $0\leq s\leq t$,
\begin{eqnarray}\label{dec.g2g3}
\begin{aligned}
\|[g_{2}(s),g_{3}(s)]\|_{L^{1}}\leq & C \|\bar{u}\|\|\nabla \bar{u}\|+\|\bar{u}\|\|B_{1}\|+\|\bar{\rho}\|\|\bar{u}\|+\|\bar{\rho}\|\|\nabla \bar{\rho}\|\\
& +C(\|\rho_{st}\nabla \bar{\rho}\|_{L^1}+\|\rho_{st} \bar{u}\|_{L^1}+\|\bar{\rho} \nabla \rho_{st}\|_{L^1} )\\
&\leq C \mathcal {E}_{N}(\bar{U}(s))+ C \delta \sqrt{\mathcal
{E}_{N}^h(\bar{U}(s))}+C \delta \|[\bar{\rho},\bar{u}]\|,
\end{aligned}
\end{eqnarray}
\begin{eqnarray}\label{dec.g2g3H}
\begin{aligned}
\|[g_{2}(s),g_{3}(s)]\|_{\dot{H}^{2}}\leq & C\mathcal
{E}_{N}(\bar{U}(s))+ C \delta \sqrt{\mathcal {E}_{N}^h(\bar{U}(s))},
\end{aligned}
\end{eqnarray}
and
\begin{eqnarray}\label{dec.g1g2}
\begin{aligned}
\|[g_{1}(s),g_{2}(s)]\|\leq & C\mathcal {E}_{N}(\bar{U}(s))+ C
\delta \sqrt{\mathcal {E}_{N}^h(\bar{U}(s))}.
\end{aligned}
\end{eqnarray}
Notice that $ \mathcal {E}_{N}(\bar{U}(s))\leq C \mathcal
{E}_{N}(\bar{V}(\sqrt{\gamma}s))$. From $\eqref{sec5.def}$ and
$\eqref{sec5.defL}$, for any $0\leq s\leq t$,
\begin{eqnarray*}
\mathcal {E}_{N}(\bar{V}(\sqrt{\gamma}s))\leq
(1+\sqrt{\gamma}s)^{-\frac{3}{2}}\mathcal
{E}_{N,\infty}(\bar{V}(\sqrt{\gamma}t)),
\end{eqnarray*}
\begin{eqnarray*}
\mathcal {E}_{N}^h(\bar{V}(\sqrt{\gamma}s))\leq
(1+\sqrt{\gamma}s)^{-\frac{5}{2}}\mathcal
{E}_{N,\infty}(\bar{V}(\sqrt{\gamma}t)),
\end{eqnarray*}
\begin{eqnarray*}
\|[\bar{\rho},\bar{u}](s)\|\leq \sqrt{L_{0}(t)}(1+s)^{-\frac{5}{4}}.
\end{eqnarray*}
Then, it follows that for $0\leq s \leq t$,
\begin{multline*}
\|[g_{2}(s),g_{3}(s)]\|_{L^{1}}\leq
C(1+\sqrt{\gamma}s)^{-\frac{3}{2}}\mathcal
{E}_{N,\infty}(\bar{V}(\sqrt{\gamma}t))\\
+C \delta (1+\sqrt{\gamma}s)^{-\frac{5}{4}}\sqrt{\mathcal
{E}_{N,\infty}(\bar{V}(\sqrt{\gamma}t))}+C\delta
\sqrt{L_{0}(t)}(1+s)^{-\frac{5}{4}},
\end{multline*}
\begin{eqnarray*}
\begin{aligned}
\|[g_{2}(s),g_{3}(s)]\|_{\dot{H}^{2}}\leq &
C(1+\sqrt{\gamma}s)^{-\frac{3}{2}}\mathcal
{E}_{N,\infty}(V(\sqrt{\gamma}t))\\
&+C \delta (1+\sqrt{\gamma}s)^{-\frac{5}{4}}\sqrt{\mathcal
{E}_{N,\infty}(\bar{V}(\sqrt{\gamma}t))},
\end{aligned}
\end{eqnarray*}
\begin{eqnarray*}
\begin{aligned}
\|[g_{1}(s),g_{2}(s)]\|\leq &
C(1+\sqrt{\gamma}s)^{-\frac{3}{2}}\mathcal
{E}_{N,\infty}(V(\sqrt{\gamma}t))\\
&+C \delta (1+\sqrt{\gamma}s)^{-\frac{5}{4}}\sqrt{\mathcal
{E}_{N,\infty}(\bar{V}(\sqrt{\gamma}t))}.
\end{aligned}
\end{eqnarray*}
Putting the above inequalities  into \eqref{sec5.decayB},
$\eqref{baruL2}$ and \eqref{rhoL2} respectively gives
\begin{multline}\label{sec5.decayB1}
 \|B_{1}(t)\|\leq C (1+t)^{-\frac{3}{4}}\Big\{
 \|[\bar{u}_{0}, E_{1,0},B_{1,0}]\|_{L^1\cap \dot{H}^{2}}+\mathcal
{E}_{N,\infty}(\bar{V}(\sqrt{\gamma}t))\\
+\delta\sqrt{L_{0}(t)}+\delta \sqrt{\mathcal
{E}_{N,\infty}(\bar{V}(\sqrt{\gamma}t))}\Big\},
\end{multline}
\begin{multline}\label{sec5.decayu}
 \|\bar{u}(t)\|\leq C (1+t)^{-\frac{5}{4}}\Big\{
 \|\bar{\rho}_{0}\|+\|[\bar{u}_{0}, E_{1,0},B_{1,0}]\|_{L^{1}\cap\dot{H}^{2}}+\mathcal
{E}_{N,\infty}(\bar{V}(\sqrt{\gamma}t))\\
+\delta\sqrt{L_{0}(t)}+\delta \sqrt{\mathcal
{E}_{N,\infty}(\bar{V}(\sqrt{\gamma}t))}\Big\},
\end{multline}
\begin{eqnarray}\label{sec5.decayrho}
\begin{aligned}
 \|\bar{\rho}(t)\|\leq C (1+t)^{-\frac{5}{4}}\Big\{
 \|[\bar{\rho}_{0},u_{0}]\|+\mathcal
{E}_{N,\infty}(\bar{V}&(\sqrt{\gamma}t))\\
&+\delta \sqrt{\mathcal
{E}_{N,\infty}(\bar{V}(\sqrt{\gamma}t))}\Big\}.
\end{aligned}
\end{eqnarray}
The definition of $L_{0}(t)$, \eqref{sec5.decayu} and
\eqref{sec5.decayrho} further imply that
\begin{multline}\label{bou.L}
L_{0}(t)\leq C\|[\bar{\rho}_{0},u_{0}]\|^{2}+C\|[\bar{u}_{0},
E_{1,0},B_{1,0}]\|_{L^{1}\cap\dot{H}^{2}}^{2}\\
+C\left[\mathcal
{E}_{N,\infty}(\bar{V}(\sqrt{\gamma}t))\right]^{2}+C\delta^{2}\mathcal
{E}_{N,\infty}(\bar{V}(\sqrt{\gamma}t)),
\end{multline}
where we have used that $\delta$ is small enough. Plugging the above
estimate into  \eqref{sec5.decayB1} implies $\eqref{lem.tildeB}$,
since $\|\bar{ B}(t)\|\leq C \| B_{1}(t/\sqrt{\gamma})\|$ and
$[\bar{\rho},\bar{u},E_{1},B_{1}]$ is equivalent with
$[\bar{\sigma},\bar{v},\bar{E},\bar{B}]$ up to a positive constant.
This completes the proof of Lemma $\ref{lem.Bsigma}$.
\end{proof}

Now, the rest is to prove the uniform-in-time bound of $\mathcal
 {E}_{N,\infty}(\bar{V}(t))$ which yields the time-decay rates of the
 Lyapunov functionals $\mathcal
 {E}_{N}(\bar{V}(t))$ and $\mathcal
 {E}_{N}^h(\bar{V}(t))$ thus $\|\bar{V}(t)\|_{N}^{2}$, $\|\nabla\bar{V}(t)\|_{N-1}^{2}$.
In fact, by taking $\ell =\frac{3}{2}+\epsilon$ in
$\eqref{sec5.EDEhD}$ with $\epsilon>0$ small enough, one has
\begin{multline*}
(1+t)^{\frac{3}{2}+\epsilon}\mathcal {E}_{N}(\bar{V}(t))+
 \int_{0}^{t}(1+s)^{\frac{3}{2}+\epsilon}\mathcal {D}_{N}(\bar{V}(s))d s\\
 +(1+t)^{\frac{5}{2}+\epsilon}\mathcal
{E}_{N}^h(\bar{V}(t))+
 \int_{0}^{t}(1+s)^{\frac{5}{2}+\epsilon}\mathcal {D}_{N}^h(\bar{V}(s))d s \\
 \leq  C \mathcal {E}_{N+3}(\bar{V}_{0})+ C \int_{0}^{t}(1+s)^{\frac{1}{2}+\epsilon}\|
  \bar{B}(s)\|^{2}d s.
\end{multline*}
Here, using $\eqref{lem.tildeB}$ and the fact $\mathcal
 {E}_{N,\infty}(\bar{V}(t))$ is non-decreasing in $t$, it further holds
 that
\begin{eqnarray*}
\begin{aligned}
 \int_{0}^{t}(1+s)^{\frac{1}{2}+\epsilon}\|
  \bar{B}(s)\|^{2}d s\leq
  C(1+t)^{\epsilon}\Big\{\|[\bar{\sigma}_{0},&\bar{v}_{0}]\|^{2}+
 \|[\bar{v}_{0},\bar{E}_{0},\bar{B}_{0}]\|^2_{L^1\cap \dot{H}^{2}}\\
 &+[\mathcal
{E}_{N,\infty}(\bar{V}(t))]^2 +\delta ^2 \mathcal
{E}_{N,\infty}(\bar{V}(t))\Big\}.
\end{aligned}
\end{eqnarray*}
Therefore, it follows that
\begin{multline*}
(1+t)^{\frac{3}{2}+\epsilon}\mathcal
{E}_{N}(\bar{V}(t))+(1+t)^{\frac{5}{2}+\epsilon}\mathcal
{E}_{N}^h(\bar{V}(t))\\
+ \int_{0}^{t}(1+s)^{\frac{3}{2}+\epsilon}\mathcal
{D}_{N}(\bar{V}(s))d s
 + \int_{0}^{t}(1+s)^{\frac{5}{2}+\epsilon}\mathcal {D}_{N}^h(\bar{V}(s))d s \\
 \leq  C \mathcal {E}_{N+3}(\bar{V}_{0})+ C (1+t)^{\epsilon}\left(\|[\bar{\sigma}_{0},\bar{v}_{0}]\|^{2}+
 \|[\bar{v}_{0},\bar{E}_{0},\bar{B}_{0}]\|^2_{L^1\cap \dot{H}^{2}}\right.\\
 \left.+[\mathcal
{E}_{N,\infty}(\bar{V}(t))]^2 +\delta^2 \mathcal
{E}_{N,\infty}(\bar{V}(t))\right),
\end{multline*}
which implies
\begin{multline*}
(1+t)^{\frac{3}{2}}\mathcal
{E}_{N}(\bar{V}(t))+(1+t)^{\frac{5}{2}}\mathcal
{E}_{N}^h(\bar{V}(t))
 \leq  C \Big\{ \mathcal {E}_{N+3}(\bar{V}_{0})+
\|[\bar{v}_{0},\bar{E}_{0},\bar{B}_{0}]\|^2_{L^1}\\
 +[\mathcal
{E}_{N,\infty}(\bar{V}(t))]^2 +\delta ^2 \mathcal
{E}_{N,\infty}(\bar{V}(t))\Big\},
\end{multline*}
and thus
\begin{eqnarray}\label{ENb}
\mathcal {E}_{N,\infty}(\bar{V}(t))
 \leq C \left( \epsilon_{N+3}(\bar{V}_{0})^{2}+
\mathcal {E}_{N,\infty}(\bar{V}(t))^{2}\right).
\end{eqnarray}
Here, we have used that $\delta$ is small enough. Recall the
definition of $\epsilon_{N+3}(\bar{V}_{0})$, since
$\epsilon_{N+3}(\bar{V}_{0})>0$ is sufficiently small, $\mathcal
{E}_{N,\infty}(\bar{V}(t)) \leq C \epsilon_{N+3}(\bar{V}_{0})^{2}$
holds true for any $t\geq 0$, which implies
\begin{eqnarray}\label{UN}
\|\bar{V}(t)\|_{N} \leq C \mathcal {E}_{N}(\bar{V}(t))^{1/2}
 \leq C  \epsilon_{N+3}(\bar{V}_{0})(1+t)^{-\frac{3}{4}},
\end{eqnarray}
\begin{eqnarray}\label{nablaUN}
\|\nabla\bar{V}(t)\|_{N-1} \leq C \mathcal
{E}_{N}^{h}(\bar{V}(t))^{1/2}
 \leq C  \epsilon_{N+3}(\bar{V}_{0})(1+t)^{-\frac{5}{4}}.
\end{eqnarray}
The definition of $L_{0}(t)$, the uniform-in-time bound of $\mathcal
{E}_{N,\infty}(\bar{V}(t))$ and $\eqref{bou.L}$ show that
\begin{eqnarray*}
\|[\bar{\rho},\bar{u}](t)\|
 \leq C  \epsilon_{N+3}(\bar{V}_{0})(1+t)^{-\frac{5}{4}}.
\end{eqnarray*}
In addition, applying the $L^{2}$ linear estimate on $E$ in
$\eqref{col.decay1}$ to the mild form $\eqref{sec5.U}$,
\begin{eqnarray*}
&&\begin{aligned} \|E_{1}(t)\|\leq C (1+t)^{-\frac{5}{4}}
 \|[\bar{u}_{0},& E_{1,0},B_{1,0}]\|_{L^1\cap \dot{H}^{3}}\\
  &+C \int_{0}^{t}(1+t-s)^{-\frac{5}{4}}\|[g_{2}(s),g_{3}(s)]\|_{L^{1}\cap
  \dot{H}^{3}}ds.
 \end{aligned}
\end{eqnarray*}
Since by $\eqref{UN}$ and $\eqref{nablaUN}$, similar to obtaining
$\eqref{dec.g2g3}$ and $\eqref{dec.g2g3H}$, we have
\begin{eqnarray*}
\begin{aligned}
&\|[g_{2}(s),g_{3}(s)]\|_{L^{1}\cap \dot{H}^{3}}\leq
C\|\bar{U}(t)\|^{2}_{4}+ C\delta\|\nabla
\bar{U}(t)\|_{3}+C\delta\|[\bar{\rho},\bar{u}]\|\leq
C\epsilon_{7}(\bar{V}_{0})(1+t)^{-\frac{5}{4}},
\end{aligned}
\end{eqnarray*}
it follows that
\begin{eqnarray}\label{uL2}
\|E_{1}(t)\| \leq C\epsilon_{7}(\bar{V}_{0})(1+t)^{-\frac{5}{4}}.
\end{eqnarray}
This completes Theorem \ref{Corolary}.

\medskip
\medskip
\noindent{\bf Acknowledgements:}\ \  The first author Qingqing Liu
would like to thank Dr. Renjun Duan for his guidance and continuous
help. The research was supported by the National Natural Science
Foundation of China $\#$11071093, the PhD specialized grant of the
Ministry of Education of China $\#$20100144110001, and the Special
Fund for Basic Scientific Research  of Central Colleges
$\#$CCNU10C01001, $\#$CCNU12C01001.

\bigbreak


\begin{thebibliography}{99}

\bibitem{Besse} C. Besse, J. Claudel, P. Degond, et al., A model hierarchy for ionospheric plasma modeling,
Math. Models Methods Appl. Sci., 14(2004), 393-415.


\bibitem{Chae} D. Chae and E. Tadmor, On the finite time blow-up of the Euler-Poisson equations
in $\mathbb{R}^{n}$, Commun. Math. Sci., 6(2008), 785-789.

\bibitem{Chen} G.Q. Chen, J.W. Jerome and D.H. Wang, Compressible Euler-Maxwell equations,
Transport Theory Statist. Phys., 29(2000), 311-331.

\bibitem{Deng} Y.B. Deng, T.P. Liu, T. Yang and Z.A. Yao, Solutions of Euler-Poisson equations for
gaseous stars, Arch. Ration. Mech. Anal., 164(2002), 261-285.


\bibitem{Duan} R.J. Duan, Global
smooth flows for the compressible Euler-Maxwell system. The
relaxation case. J. Hyperbolic Differ. Equ., 8(2011), 375-413.

\bibitem{DLUY} R.J. Duan, H.X. Liu, S. Ukai and T. Yang, Optimal $L^p$-$L^q$ convergence rates
for the compressible Navier-Stokes equations with potential force,
J. Differential Equations, 238(2007), 220-233.

\bibitem{DLZ} R.J. Duan, Q.Q. Liu, C.J. Zhu, The Cauchy problem on the compressible two-fluids Euler-Maxwell equations,
SIAM J. Math. Anal.,  44(2012), 102-133.

\bibitem{DS} R.J. Duan and R.M. Strain, Optimal time decay of the
Vlasov-Poisson-Boltzmann system in $\mathbb{R}^{3}$, Arch. Ration.
Mech. Anal., 199(2011), 291-328.


\bibitem{RSY} R.J. Duan, S. Ukai and T. Yang, A combination of energy method and spectral
analysis for study of equations of gas motion, Front. Math. China,
4(2009), 253-282.


\bibitem{RSYZ} R.J. Duan, S. Ukai, T. Yang and H.J Zhao, Optimal convergence rates for
the compressible Navier-Stokes equations with potential forces,
Math. Models Methods Appl. Sci., 17(2007), 737-758.

\bibitem{RY} R.J. Duan and T. Yang, Stability of the one-species Vlasov-Poisson-Boltzmann
system, SIAM J. Math. Anal., 41(2009/10), 2353-2387.






 \bibitem{Guo}  Y. Guo, Smooth irrotational flows in the large to the Euler-Poisson system in $\mathbb{R}^{3+1}$,
Comm. Math. Phys., 195(1998), 249-265.

\bibitem{GuoPausader} Y. Guo, B. Pausader, Global smooth ion dynamics in the Euler-Poisson
system, Comm. Math. Phys., 303(2011), 89-125.

\bibitem{HP} M.L. Hajjej, Y.J. Peng, Initial layers and zero-relaxation limits of Euler-Maxwell
equations, J. Differential Equations, 252(2012), 1441-1465.



\bibitem{Jerome} J.W. Jerome, The Cauchy problem for compressible hydrodynamic-Maxwell
systems: a local theory for smooth solutions, Differential Integral
Equations, 16(2003), 1345-1368.





 \bibitem{Kato} T. Kato, The Cauchy problem for quasi-linear symmetric hyperbolic systems, Arch.
Rational Mech. Anal., 58(1975), 181-205.




 \bibitem{Luo}T. Luo, R. Natalini and Z.P. Xin, Large time behavior of the solutions to a
 hydrodynamic model for semiconductors, SIAM J. Appl. Math., 59(1999), 810-830.


\bibitem{Smoller}T. Luo and J. Smoller, Existence and non-linear stability
of rotating star solutions of the compressible Euler-Poisson
equations, Arch. Ration. Mech. Anal., 191(2009), 447-496.




\bibitem{Peng }Y.J. Peng and S. Wang, Convergence of compressible Euler-Maxwell equations to
incompressible Euler equations, Comm. Partial Differential
Equations, 33(2008), 349-376.


\bibitem{PW1} Y.J. Peng, S. Wang, Convergence of compressible Euler-Maxwell
equations to compressible Euler-Poisson equations, Chin. Ann. Math.
Ser. B, 28(2007), 583-602.

\bibitem{PW2}Y.J. Peng, S. Wang, Rigorous
derivation of incompressible e-MHD equations from compressible
Euler-Maxwell equations, SIAM J. Math. Anal., 40(2008) 540-565.


\bibitem{Rishbeth} H. Rishbeth and O.K. Garriott,
Introduction to Ionospheric Physics, Academic Press, 1969.





 \bibitem{Stein} E.M. Stein, Singular Integrals and Differentiability Properties of
Functions, Princeton Mathematical Series, No. 30 Princeton
University Press, Princeton, N.J. 1970 xiv+290 pp.

 \bibitem{Taylor} M.E. Taylor, Partial Differential Equations, I. Basic Theory, Springer, New York,
1996.

\bibitem{UK}
Y. Ueda and S. Kawashima, {Decay property of regularity-loss type
for the Euler-Maxwell system}, Methods Appl. Anal., 18(2011),
245-268.

\bibitem{USK}
Y. Ueda, S. Wang and S. Kawashima, {Dissipative Structure of the
Regularity-Loss Type and Time Asymptotic Decay of Solutions for the
Euler-Maxwell System}, SIAM J. Math. Anal., 44(2012), 2002-2017.





\end{thebibliography}
\end{document}